\documentclass[a4paper,11pt]{amsart}
\usepackage[frenchb]{babel}
\usepackage{amssymb,xypic,latexsym,graphicx,caption2}
\usepackage[matrix,arrow]{xy}
\newtheorem{theoreme}{Th\'eor\`eme}
\newtheorem{cor}{Corollaire} 
\newtheorem{lemme}{Lemme}

\newtheorem{prop}{Proposition}

\newcommand{\cucu}{\hfill $\Box$}
\newcommand{\SSS}{{{\mathcal N}^{d,n}}}
\newcommand{\YYY}{{\mathcal V}^{d}}

\newcommand{\HH}{{\rm H}}
\newcommand{\hh}{{\rm h}}

\newcommand{\ZZ}{{\mathbb Z}}
\newcommand{\QQ}{{\mathbb Q}}
\newcommand{\NN}{{\mathbb N}}
\newcommand{\RR}{{\mathbb R}}
\newcommand{\C}{{\mathbb C}}
\newcommand{\PP}{{\mathbb P}}

\newcommand{\codim}{\text{codim}\,}
\newcommand{\kker}{\text{Ker}\,}

\newcommand{\anneau}{{\mathcal O}}
\newcommand{\preuve}{\noindent{\it Preuve. --- \ }}

\newcommand{\ev}{{\rm ev}}

\newcommand{\tube}{{\rm tube}}

\begin{document}
\def\figurename{{\sc Dessin}}
\title{Sur les vari\'et\'es de Hodge des hypersurfaces}
\author{Ania Otwinowska}
\date{}
\maketitle


\addtocounter{section}{-1}

\section{Introduction}
 
\subsection{Avant-propos}
Un th\'eor\`eme c\'el\`ebre de Noether affirme que le groupe de Picard 
d'une surface g\'en\'erale de $\PP^3_{\C}$ de degr\'e sup\'erieur ou
\'egal \`a $4$ est 
engendr\'e par la classe d'une section hyperplane. Identifions le groupe
de Picard d'une surface lisse $X\subset\PP^3_{\C}$ \`a l'espace
$\HH^{1,1}(X)\cap\HH^2(X,\ZZ)$ des classes de Hodge~; le th\'eor\`eme de
Noether admet alors la g\'en\'eralisation naturelle suivante. Soit $Y$ une
vari\'et\'e complexe projective lisse de dimension $2n+1$ et ${\mathcal L}$ un
faisceau inversible suffisamment ample sur $Y$. Alors une hypersurface
g\'en\'erale $X$ de $|{\mathcal L}|=\PP\HH^0(Y,{\mathcal L})$ ne 
poss\`ede pas de classe de Hodge non triviale~: 
nous avons $\HH^{n,n}(X)_{\ev}\cap\HH^{2n}(X,\QQ)=0$, o\`u
$\HH^{n,n}(X)_{\ev}=\HH^{n,n}(X)\cap\HH^{2n}(X,\C)_{\ev}$ d\'esigne un 
suppl\'ementaire de l'image de $\HH^{n,n}(Y)$ dans $\HH^{n,n}(X)$ 
(la d\'efinition pr\'ecise est donn\'ee \`a la section~\ref{introtopo}).   

La motivation principale de cet article est l'\'etude du {\em lieu de
Noether-Lefschetz} param\'etrant les hypersurfaces lisses 
$X\in |{\mathcal L}|$ qui
poss\`edent une classe de Hodge non triviale~; c'est une r\'eunion
d\'enombrable de sous-vari\'et\'es alg\'ebriques de l'ouvert 
${\mathcal V}_{\mathcal L}\subset|{\mathcal L}|$ param\'etrant 
les hypersurfaces lisses. Il
refl\`ete le comportement de la cohomologie rationnelle $\HH^{2n}(X,\QQ)$
vis-\`a-vis de la variation des structures de Hodge sur la cohomologie
\'evanescente \`a coefficients complexes $\HH^{2n}(X,\C)_{\ev}$. Son
importance tient en grande partie \`a la conjecture de Hodge qui affirme
que l'espace des classes de Hodge de $X$ est engendr\'e par les classes
des sous-vari\'et\'es alg\'ebriques de $X$~; les composantes du lieu de
Noether-Lefschetz devraient donc param\'etrer des familles non triviales
de sous-vari\'et\'es alg\'ebriques des hypersurfaces de 
${\mathcal V}_{\mathcal L}$.

Soit plus g\'en\'eralement $Y$ une vari\'et\'e complexe projective 
lisse de dimension $N+1$, ${\mathcal L}$ un faisceau inversible 
suffisamment ample sur $Y$, $X\in |{\mathcal L}|$ une hypersurface 
lisse et $\lambda\in F^k\HH^{N}(X,\C)_{\ev}$ une classe de cohomologie
\'evanescente, o\`u $F^{\bullet}$ d\'esigne la filtration de Hodge et
$k\in\{1,\dots,N\}$. On appelle {\em vari\'et\'e de Hodge} associ\'ee \`a
$\lambda$ dans un voisinage contractile $U\subset |{\mathcal L}|$ de $X$ le
ferm\'e analytique ${\mathcal N}_{\lambda,U}^{{\mathcal L},k}$ o\`u 
l'image par transport plat de la classe $\lambda$ reste dans
$F^k\HH^{N}(X,\C)_{\ev}$. Dans le cas $\lambda\in\HH^{N}(X,\QQ)$ et
$2k=N$, la classe $\lambda$ est de Hodge, et le lieu 
${\mathcal N}_{\lambda,U}^{{\mathcal L},k}$ est inclus dans le lieu de
Noether-Lefschetz.   

Nous montrons que {\em si ${\mathcal L}$ est suffisamment ample et si
la codimension de la vari\'et\'e de Hodge 
${\mathcal N}_{\lambda,U}^{{\mathcal L},k}$ associ\'ee \`a une 
classe  $\lambda\in F^k\HH^{N}(X,\C)_{\ev}$ est suffisamment petite,
la classe $\lambda$ est une combinaison lin\'eaire \`a coefficients
complexes des projet\'es sur la cohomologie \'evanescente des classes
de sous-vari\'et\'es alg\'ebriques de $X$ de petit degr\'e}. En
particulier, la petite codimension de 
${\mathcal N}_{\lambda,U}^{{\mathcal L},k}$ force la
dimension de $X$ \`a \^etre paire~: nous avons $N=2k$. 

Comme corollaire, nous obtenons que {\em les composantes de plus
petites codimensions du lieu de Noether-Lefschetz sont d\'efinies par
des classes des sous-vari\'et\'es alg\'ebriques}, comme le pr\'edit la
conjecture de Hodge. 
Cependant les classes de sous-vari\'et\'es alg\'ebriques 
(de petit degr\'e) sont ici caract\'eris\'ees par leur seul comportement
vis-\`a-vis de la variation des structures de Hodge sur ${\mathcal
V}_{\mathcal L}$, sans qu'intervienne leur caract\`ere rationnel.



\subsection{\'Enonc\'e du th\'eor\`eme principal}
Soit $Y$ une vari\'et\'e complexe projective lisse de dimension $N+1$
munie d'un faisceau inversible tr\`es ample $\anneau (1)$, soit
$d\in\NN^*$ et soit
$\YYY\subset \HH^0(Y,\anneau (d))$ l'ouvert param\'etrant (\`a la
multiplication par un scalaire pr\`es) les hypersurfaces projectives
lisses de $Y$ de classe $c_1(\anneau (d))$. Pour tout $F\in\YYY$ nous
notons $j_F:X_F\to Y$ l'immersion ferm\'ee de l'hypersurface
correspondante. 

\subsubsection{Cohomologie \'evanescente et sous-vari\'et\'es de $X_F$} 
\label{introtopo}
Nous notons
  $$\HH^{N}(X_F,\C)_{\ev}=
  \kker(g:\HH^{N}(X_F,\C)\rightarrow\HH^{N+2}(Y,\C)),$$ 
la cohomologie \'evanescente de $X_F$, o\`u $g$ d\'esigne le morphisme de 
Gysin. Le th\'eor\`eme de Lefschetz difficile affirme que nous avons 
une somme directe orthogonale pour le cup-produit
$$\HH^{N}(X_F,\C)=\HH^{N}(X_F,\C)_{\ev}\oplus j_F^*\HH^{N}(Y,\C).$$ 

Si $W$ est une sous-vari\'et\'e de $X_F$ (c'est-\`a-dire un sous-sch\'ema 
ferm\'e r\'eduit) de dimension pure $n$, chaque composante irr\'eductible 
de $W$ d\'efinit par dualit\'e de Poincar\'e une classe dans 
$\HH^{2N-2n}(X_F,\C)$. Lorsque $2n=N$, nous notons $H(W)$ le sous-espace 
vectoriel de $\HH^{N}(X_F,\C)_{\ev}$ engendr\'e par les projections 
orthogonales de ces classes.

\subsubsection{Vari\'et\'es de Hodge}\label{introNL}    
Nous avons une famille de structures de Hodge  
(${\mathcal H}^{N}_{\ZZ} =R^{N}\pi_*\ZZ$~; 
${\mathcal H}^{N}={\mathcal H}^{N}_{\ZZ}\otimes\anneau_{\YYY}$~;
$F^{\bullet}{\mathcal H}^{N}$) sur $\YYY$, telle que pour tout $F\in\YYY$
l'espace $\HH^{N}(X_F,\C)$ s'identifie canoniquement \`a la fibre en 
$F$ du faisceau localement constant ${\mathcal H}^{N}_{\ZZ}\otimes
\C$. Pour tout voisinage ouvert simplement connexe $U\subset \YYY$ de
$F$, le faisceau ${\mathcal H}^{N}_{\ZZ}\otimes \C$ restreint \`a $U$
est constant et toute classe $\lambda_F\in\HH^{N}(X_F,\C)$ d\'efinit une
section $\lambda \in ({\mathcal H}^{N}_{\ZZ}\otimes \C)(U)$. 
Pour tout $G\in U$ nous notons $\lambda_G\in\HH^{N}(X_G,\C)$ 
la valeur de la section $\lambda$ en $G$~: c'est
l'image de la classe $\lambda_F$ par transport plat dans $U$.

Soit $F\in\YYY$, soit $n\in\{0,\dots,N\}$ et soit
$\lambda_F\in F^{n}\HH^{N}(X_F,\C)_{\ev}$ une classe de cohomologie
non nulle. Nous introduisons le ferm\'e analytique
$${\mathcal N}_{\lambda,U}^{d,n}=\{G\in U \ | \ 
\lambda_G\in F^n\HH^{N}(X_G,\C)\}.$$

\subsubsection{Cas N=2n}
Les vari\'et\'es ${\mathcal N}_{\lambda,U}^{d,n}$ jouent un r\^ole
particuli\`erement important lorsque $2n=N$. On appelle alors lieu de
Noether-Lefschetz, que nous notons $\SSS$, la r\'eunion des vari\'et\'es
${\mathcal N}_{\lambda,U}^{d,n}$ pour tous les couples $(\lambda,U)$ tels
que la classe  $\lambda\in ({\mathcal H}^{2n}_{\ZZ}\otimes \C)(U)$ est
rationnelle, c'est-\`a-dire appartient \`a 
$({\mathcal H}^{2n}_{\ZZ}\otimes\QQ)(U)$~; nous avons encore
$$\SSS=\{G\in\YYY\ |\ 
\HH ^{2n}(X_G,\QQ)\cap F^n\HH^{2n}(X_G,\C)_{\ev}\not=0\}.$$  
Le lieu  
${\mathcal N}_{\lambda}^{d,n}=\bigcup_U {\mathcal N}_{\lambda,U}^{d,n}$ 
d\'efini par recollement est une sous-vari\'et\'e {\em alg\'ebrique} de
$\YYY$ et $\SSS$ est une r\'eunion d\'enombrable de vari\'et\'es 
{\em alg\'ebriques} de $\YYY$ ({\it cf.}~\cite{CDK}). 
La conjecture de Hodge pr\'edit que toute classe de 
$\HH ^{2n}(X_F,\QQ)\cap F^n\HH^{2n}(X_F,\C)$ est
alg\'ebrique, c'est-\`a-dire combinaison lin\'eaire \`a coefficients
rationnels de classes de sous-vari\'et\'es alg\'ebriques de $X_F$. 
Si nous supposons la conjecture de Hodge vraie pour $Y$, 
le lieu de Noether-Lefschetz param\`etre donc les hypersurfaces 
pour lesquelles la conjecture de Hodge est non triviale.

\subsubsection{\'Enonc\'e} Nous montrons ici~:

\begin{theoreme}\label{nlh+r}
Soient $n\in\{1,\dots,[N/2]\}$ et $b\in \NN^*$ des entiers,  
$F\in \YYY$, $U$ un voisinage ouvert simplement connexe de $F$ dans
$\YYY$ et $\lambda_F\in F^n\HH^{N}(X_F,\C)_{\ev}$ une classe non
nulle v\'erifiant 
$${\rm codim}\,({\mathcal N}_{\lambda,U}^{d,n},U) \leq b\frac{d^n}{n!}.$$ 
Pour tout $\varepsilon \in \RR^*_+$, il existe une constante
$C\in\RR^*_+$ qui ne d\'epend que de $Y$, de $\anneau(1)$ et de
$\varepsilon$, telle que si $d\geq Cb^{a}$, o\`u
$a=2^{\hh^0(Y,\anneau(1))-1}$, alors pour tout
$n$, $F$, $U$ et $\lambda_F$ comme ci-dessus 
\begin{enumerate}
\item $2n=N$ et
\item il existe une sous-vari\'et\'e $W$ de $X_F$ de
dimension $n$ et de degr\'e inf\'erieur  \`a
$(1+\varepsilon)b$, telle que $\lambda\in H(W)$.  
\end{enumerate}
\end{theoreme}

\noindent
Nous pensons que le th\'eor\`eme~\ref{nlh+r} reste vrai pour 
$n\in\{[N/2]+1,\dots,N\}$.

\medskip 

Le cas $2n=N$ du th\'eor\`eme~\ref{nlh+r} implique imm\'ediatement le 
corollaire suivant, qui est une version asymptotique tr\`es faible de la
conjecture de Hodge pour les hypersurfaces.     

\begin{cor}
Nous supposons $Y$ de dimension $2n+1$ et que toute classe dans $\HH
^{2n}(Y,\QQ)\cap F^n\HH^{2n}(Y,\C)$ est alg\'ebrique. Alors il existe 
une constante $C\in\RR_+^*$ telle que toute classe 
$\lambda_F\in\HH ^{2n}(X_F,\QQ)\cap F^n\HH^{2n}(X_F,\C)$ v\'erifiant
$${\rm codim}\,({\mathcal N}_{\lambda}^{d,n},\YYY)
\leq\left(\frac{d}{C}\right)^\frac{1}{a}\frac{d^n}{n!},$$
o\`u $a=2^{\hh^0(Y,\anneau(1))}$, est alg\'ebrique.
\end{cor}

L'hypoth\`ese de rationnalit\'e n'est utilis\'ee ici que pour la 
composante non \'evanescente de la classe $\lambda_F$. Plut\^ot que de
plaider en faveur de la conjecture de Hodge, le th\'eor\`eme~\ref{nlh+r}
met donc en \'evidence le comportement exceptionnel des classes
alg\'ebriques vis-\`a-vis de la variation des structures de Hodge.

\medskip

Le th\'eor\`eme~\ref{nlh+r} est l'aboutissement de l'\'etude des composantes
des petite codimension du lieu de Noether-Lefschetz men\'ee dans
\cite{NLgreen}, \cite{NLclaire}, \cite{NLclaire2} et~\cite{NL}. Il 
compl\`ete et g\'en\'eralise~\cite{NL}, o\`u nous
avons montr\'e dans le cas $Y=\PP_{\C}^{2n+1}$ que sous les hypoth\`eses
du th\'eor\`eme~\ref{nlh+r} les vari\'et\'es ${\mathcal N}_{\lambda,U}^{d,n}$
sont constitu\'ees par des hypersurfaces contenant une sous-vari\'et\'e $W$ 
de dimension $n$ et de degr\'e inf\'erieur  \`a $(1+\varepsilon)b$, 
mais o\`u nous ne savions pas montrer $\lambda\in H(W)$.
Le cas $Y=\PP_{\C}^{N+1}$, $n=[N/2]$ du th\'eor\`eme~\ref{nlh+r} r\'esout
asymptotiquement la conjecture 1 de~\cite{NL}.

\subsection{Structure de la preuve}

\subsubsection{Invariants alg\'ebriques des vari\'et\'es de Hodge}
La th\'eorie de Griffiths~\cite{grI} et de Carlson, Green, Griffiths 
et Harris~\cite{ivhs} ram\`ene l'\'etude des variations 
infinit\'esimales de structures de Hodge d'hypersurfaces suffisamment
amples \`a un probl\`eme alg\'ebrique (section~\ref{rappels}). 
Elle d\'ecrit notamment le lieu ${\mathcal N}_{\lambda,U}^{d,n}$ 
\`a l'ordre~1 au voisinage de $F$, et cette description est exploit\'ee 
dans~\cite{ivhs2},~\cite{NLgreen},~\cite{NLclaire} et~\cite{NLclaire2}.

En raffinant cette \'etude nous d\'ecrivons le lieu 
${\mathcal N}_{\lambda,U}^{d,n}$ au voisinage de $F$ par une suite 
$E_r$ d'id\'eaux gradu\'es de l'alg\`ebre 
$\bigoplus_{i\in\NN}\HH^0(Y,\anneau_Y(i))$  
(section~\ref{defo}).

\subsubsection{\'Etude locale des vari\'et\'es de Hodge}
Cette section g\'en\'eralise~\cite{NL}.

Nous traduisons d'abord les propri\'et\'es g\'eom\'etriques du lieu
${\mathcal N}_{\lambda,U}^{d,n}$ par des 
propri\'et\'es alg\'ebriques des id\'eaux $E_r$ (proposition~\ref{3prop}).

D'autre part,
l'\'etude asymptotique de la fonction de Hilbert de ces id\'eaux 
(lorsque le degr\'e $d$ de l'hypersurface tend vers l'infini) permet 
de montrer qu'il existe un entier $i$ petit devant $d$ tel 
que le lieu-base des \'el\'ements homog\`enes de degr\'e $i$ de $E_0$
est de dimension $n$ (lemme~\ref{llh}). En affinant notre \'etude,
nous montrons ensuite par des techniques purement alg\'ebriques 
que l'id\'eal $E_0$ contient l'id\'eal d'un sous-sch\'ema $V$ de $Y$ 
de dimension pure $n$ et de petit degr\'e (proposition~\ref{alh}). 

Nous supposons $d$ assez grand et nous montrons gr\^ace aux  
propositions~\ref{3prop} et~\ref{alh} que pour tout $r\in\NN^*$ 
l'id\'eal $E_r$ contient le produit d'un id\'eal fixe 
(trivial lorsque $V$ est r\'eduit) et de la puissance $(r+1)$-i\`eme de 
l'id\'eal $I_W$ du sch\'ema r\'eduit $W$ 
associ\'e \`a $V$ (proposition~\ref{I^2am}). Enfin, nous montrons 
$W\subset X_F$ (proposition~\ref{Zred}).

\subsubsection{\'Etude globale des vari\'et\'es de Hodge}
Nous \'etudions l'image par transport plat de la classe 
$\lambda_F$ dans le lieu ${\mathcal V}^d(W)$ des hypersurfaces 
lisses contenant $W$. 
Nous identifions l'espace tangent en $F$ \`a ${\mathcal V}^d(W)$ \`a 
l'espace des \'el\'ements homog\`enes de degr\'e $d$ de $I_W$.
Nous montrons ensuite gr\^ace aux propositions~\ref{I^2am} 
et~\ref{Zred} que les d\'eriv\'ees \`a tout ordre en $F$ de certaines
\'equations de ${\mathcal N}_{\lambda,U}^{d,n}$
le long du lieu ${\mathcal V}^d(W)$ sont nulles
(et m\^eme de toutes les \'equations dans le cas o\`u $V$ est r\'eduit).
Nous en d\'eduisons par un argument d'unicit\'e du prolongement 
analytique que la repr\'esentation de monodromie de
$\pi_1({\mathcal V}^d(W),F)$ engendr\'ee par $\lambda$ est contenue
dans un petit sous-espace de $H^N(X_F,\C)_{\ev}$ 
(proposition~\ref{resume}).

La fin de la preuve repose alors sur un \'enonc\'e de topologie
que nous avons montr\'e dans~\cite{topo}~: pour tout sous-sch\'ema 
ferm\'e $W$ de $X_F$ dont l'id\'eal est engendr\'e en degr\'e 
petit devant $d$, l'action de monodromie de 
$\pi_1({\mathcal V}^d(W),F)$ est triviale sur $H(W)$ et  
irr\'eductible sur son orthogonal dans $H^N(X_F,\C)_{\ev}$ 
(proposition~\ref{mono}). 
La proposition~\ref{resume} implique alors que la projection 
orthogonale de $\lambda$ sur $H(W)^{\perp}$ est nulle,
ce qui implique le th\'eor\`eme~\ref{nlh+r}.

\subsubsection{Remarque} Dans~\cite{NL} nous conjecturons le r\'esultat du
th\'eor\`eme~\ref{nlh+r} avec une borne sur $d$ lin\'eaire en $b$, alors que  
la borne que nous obtenons n'est que polynomiale en $b$~:
l'obstruction essentielle pour obtenir une borne   
lin\'eaire vient de ce que nous ne savons pas borner lin\'eairement 
en $b$ le degr\'e du plus grand g\'en\'erateur de $I_W$~; la borne
polynomiale est obtenue gr\^ace au lemme~\ref{bamu}.

\section{Invariants alg\'ebriques des vari\'et\'es de Hodge}\label{delo}

Les notations sont celles de la section~\ref{introNL}.
Posons
$$S^i=\HH^0(Y,\anneau_Y(i)),\quad S=\bigoplus_{i\in\NN}S^i
 \qquad{\rm et}\qquad 
K^i=\HH^0(Y,\Omega^{N+1}_Y\otimes\anneau_Y(i)),\quad
K=\bigoplus_{i\in\ZZ}K^i,$$
si bien que $K$ est un $S$-module gradu\'e de type fini.
Plus g\'en\'eralement, pour tout $S$-module gradu\'e $M$ nous notons 
$M^i\subset M$ les \'el\'ements homog\`enes de degr\'e $i$.

\subsection{Rappels sur les lieux de Hodge}\label{rappels}
Nous rappelons bri\`evement une th\'eorie expos\'ee en d\'etail 
dans~\cite{grI} et~\cite{ivhs}. 

\subsubsection{Th\'eorie de Griffiths}
Nous avons une application r\'esidu surjective
$${\rm res}:\HH^{N+1}(Y\setminus X_F,\C)\to\HH^N(X_F,\C)_{\ev}.$$
Pour tout $k\in\ZZ$ et $\omega\in K^{kd}$, la forme diff\'erentielle
$\frac{\omega}{F^k}\in\HH^0(Y\setminus X_F,\Omega_{Y\setminus X_F}^{N+1})$
d\'efinit une classe 
$\left[\frac{\omega}{F^k}\right]\in\HH^{N+1}(Y\setminus X_F,\C)$.
Nous notons
$\phi_{F,k}(\omega)={\rm res}\left[\frac{\omega}{F^k}\right]$ 
son image dans $\HH^N(X_F,\C)_{\ev}$.

            Griffiths a montr\'e que si $d$ est plus grand qu'une 
constante $C\in\RR_+^*$ qui ne d\'epend que de $Y$ et de $\anneau(1)$, 
cette construction donne pour tout $k\in\NN^*$ des applications surjectives 
$$\phi_{F,k}\colon K^{kd}\to F^{N+1-k}\HH^N(X_F,\C)_{\ev}$$
o\`u nous avons pos\'e 
$F^k\HH^N(X_F,\C)_{\ev}=F^k\HH^N(X_F,\C)\cap\HH^N(X_F,\C)_{\ev}$ pour
$k\in\{1,\dots,N+1\}$ et $F^k\HH^N(X_F,\C)_{\ev}=\HH^N(X_F,\C)_{\ev}$ pour
$k\leq 0$.  

Nous supposerons dor\'enavant $d\geq C$. 

\subsubsection{\'Equations locales des vari\'et\'es de Hodge}\label{equa}
Comme $F^k\HH^N(X_F,\C)_{\ev}$ est l'orthogonal de 
$F^{N+1-k}\HH^N(X_F,\C)_{\ev}$ pour le cup-produit, la surjectivit\'e de
$\phi_{F,k}$ implique pour tout $\mu_F\in\HH^N(X_F,\C)_{\ev}$
$$\mu_F\in F^k\HH^N(X_F,\C)_{\ev}\quad\Longleftrightarrow\quad
\forall\omega\in K^{kd},\quad\mu_F\smile\phi_{F,k}(\omega)=0,$$
o\`u $\smile$ d\'esigne le cup-produit.

Supposons $\mu_F\in F^k\HH^N(X_F,\C)_{\ev}$. Alors pour tout $G\in U$
$$G\in {\mathcal N}_{\mu,U}^{d,k}\quad\Longleftrightarrow\quad
  \forall\omega\in K^{kd},\qquad\mu_G\smile\phi_{G,k}(\omega)=0.$$ 
La vari\'et\'e ${\mathcal N}_{\mu,U}^{d,k}\subset U$
est donc d\'efinie par $\dim\,(K^{kd})$ \'equations posssiblement 
redondantes.


\subsection{D\'eformations \`a tout ordre des lieux de Hodge}\label{defo}
Nous fixons d\'esormais $n\in\{1,\dots,N\}$, $F\in \YYY$, $U$ un
voisinage ouvert simplement connexe de $F$ dans $\YYY$  et
$\lambda_F \in F^n\HH^{N}(X_F,\C)_{\ev}$. 

\subsubsection{D\'efinition} Pour tout $r\in\NN$ nous posons
$$E_{r}=\left\{\,
  P\in S^{\bullet}\ |\ \forall\omega\in K^{(n+r+1)d-\bullet}, \
  \lambda_F\smile\phi_{F,n+r+1}(P\omega)=0\right\}.$$
L'id\'eal $E_0$ a \'et\'e \'etudi\'e dans~\cite{NLgreen} et l'id\'eal 
$E_1$ dans~\cite{NL}.

\subsubsection{Description alternative des id\'eaux $E_r$}\label{alte}
Nous choisissons des cycles  
$\gamma_{F,i}$, $i\in\{1,\dots, t\}$, $t=\hh_{N}(X_F)$, 
dont les classes $[\gamma_{F,i}]$ forment une base de $\HH_{N}(X_F,\QQ)$.
Nous notons $\lambda_F^{\vee}\in\HH_{N}(X_F,\C)$ l'image de 
$\lambda_F\in\HH^{N}(X_F,\C)_{\ev}$ par la dualit\'e de
Poincar\'e et $(\lambda_1,\dots,\lambda_t)$ les coordonn\'ees 
complexes de $\lambda_F^{\vee}$~: nous avons  
$\lambda_F^{\vee}=\sum_{i=1}^t\lambda_i[\gamma_{F,i}]$.

Nous choisissons ensuite des cycles ${\rm Tub}\,\gamma _i$, 
$i\in\{1,\dots, t\}$, \`a support dans 
$Y\setminus X_F$ dont les classes
$[{\rm Tub}\,\gamma _i]$ sont les images des
$[\gamma_{F,i}]$ par l'application 
$\tube:\HH_{N}(X_F,\ZZ)\to\HH_{N+1}(Y\setminus X_F,\ZZ)$.

Quitte \`a restreindre le voisinage $U\subset\YYY$ de $F$, nous
supposons que pour tout $G\in U$ et pour tout $i\in\{1,\dots, t\}$ 
nous avons $X_G\cap {\rm Tub}\,\gamma_i=\emptyset$.
La classe de ${\rm Tub}\,\gamma_i$ dans 
$\HH_{N+1}(Y\setminus X_G,\ZZ)$ s'identifie alors 
\`a l'image par l'application tube de 
$[\gamma_{G,i}]\in\HH_{N}(X_G,\ZZ)$, image de 
$[\gamma_{F,i}]\in\HH_{N}(X_F,\ZZ)$ par transport plat dans $U$. 

Pour tout $k\in\NN^*$ et pour tout $\omega\in K^{kd}$ 
nous avons donc
$$\lambda_G\smile\phi_{G,k}(\omega)=
\sum_{i=1}^t\lambda_i\int_{\gamma_{G,i}}\phi_{G,k}(\omega)=
\sum_{i=1}^t\lambda_i\int_{{\rm Tub}\,\gamma_i}\frac{\omega}{G^{k}}.$$ 
En particulier, nous en d\'eduisons pour tout $r\in\NN$
$$E_{r}=\left\{\,
  P\in S^{\bullet}\ |\ \forall\omega\in K^{(n+r+1)d-\bullet}, \
  \sum_{i=1}^t\lambda_i\int_{{\rm Tub}\,\gamma_i}\frac{P\omega}{F^{n+r+1}}=0
  \right\}.$$

\subsubsection{Interpr\'etation g\'eom\'etrique}
L'id\'eal $E_r\subset S$ d\'ecrit le lieu
${\mathcal N}_{\lambda,U}^{d,n}$ \`a l'ordre $r+1$ au voisinage de
$F$. Plus pr\'ecis\'ement, si 
$\tau_F\colon T_F U\to S^d$ d\'esigne l'isomorphisme canonique, pour tout 
$(\vec{v}_1,\dots,\vec{v}_{r+1})\in(T_FU)^{r+1}$
les assertions suivantes sont \'equivalentes~:
\begin{enumerate}
\item 
$\forall\omega\in K^{nd},\ \frac{\partial\ }{\partial \vec{v}_1}\cdots
\frac{\partial\ }{\partial \vec{v}_{r+1}}
(G\mapsto\lambda_G\smile\phi_{G,n}(\omega))(F)=0$
\item 
$\tau_F(\vec{v}_1)\cdots\tau_F(\vec{v}_{r+1})\in E_{r+1}^{d(r+1)}$.
\end{enumerate}
Remarquons que $E_r\neq S$ d\`es que 
$\lambda_F \not\in F^{n+r+1}\HH^{N}(X_F,\C)_{\ev}$.

\subsubsection{Point de vue dual} 
Pour tout $r\in\NN$ nous pouvons d\'efinir le module gradu\'e 
$$M_r=\left\{\,\omega\in K^{\bullet}\ |\ \forall
  P\in S^{(n+r+1)d-\bullet}, \ 
\lambda_F\smile\phi_{F,n+r+1}(P\omega)=0\right\}.$$
L'espace $\phi_{F,n+r+1}\left(M_r^{(n+r+1)d}\right)$ est alors 
l'espace $\bigoplus_{i\in\{n+r+1,\dots,N\}}H^{i,N-i}(\lambda_F)$
\'etudi\'e dans~\cite{ivhs} et~\cite{ivhs2}. 

L'id\'eal $E_r$ est l'annulateur du module $K/M_r$ et la forme 
lin\'eaire $\omega\mapsto\lambda_F\smile\phi_{F,n+r+1}(\omega)$
sur $K^{(n+r+1)d}$ induit une dualit\'e parfaite 
$$(K/M_r)^{\bullet}\otimes(S/E_r)^{(n+r+1)d-\bullet}\to\C.$$

\subsubsection{Les id\'eaux $E_{r,i}$} \label{defogo}
Le $S$-module gradu\'e $K$ est de type fini~; il est donc 
engendr\'e en degr\'e inf\'erieur  \`a un entier $m\in\NN$ 
qui ne d\'epend que de $Y$ et de $\anneau_Y(1)$.

Posons $u=\dim K^m$ et choisissons une base 
$(\omega_1,\dots,\omega_u)$ de $K^m$. 
Pour tout $i\in\{1,\dots,u\}$ posons
$$E_{r,i}=\left\{ P\in S^{\bullet}\ |\ \forall Q\in 
  S^{(n+r+1)d-m-\bullet}, \, 
  \lambda_F\smile\phi_{F,n+r+1}(PQ\omega_i)=0\right\}.$$
L'id\'eal $E_{r,i}$ est l'annulateur du module 
$S\omega_i/(M_r\cap S\omega_i)$. 
Lorsque $E_{r,i}\neq S$, la forme lin\'eaire 
$P\mapsto\lambda_F\smile\phi_{F,n+r+1}(P\omega_i)$ sur 
$S^{(n+r+1)d-m}$ est non nulle et l'id\'eal $E_{r,i}$ est 
Gorenstein de degr\'e du socle $(n+r+1)d-m$ (voir section~\ref{gore} 
pour la d\'efinition d'un id\'eal Gorenstein). Nous avons
$$E_r=\bigcap_{i=1}^uE_{r,i}.$$

\section{Description locale des vari\'et\'es de Hodge}

\subsection{Propri\'et\'es des id\'eaux $E_r$ issues de la g\'eom\'etrie}
\label{geo}
Dans cette section nous \'enon\c cons 
les propri\'et\'es alg\'ebriques des id\'eaux $E_r$ qui r\'esultent 
de leur interpr\'etation g\'eom\'etrique. 

Les notations sont celles de la section~\ref{defo}. 
Pour tout id\'eal $I\subset S$ et pour $r\in\NN$ nous notons $(I)^r$
la puissance r-i\`eme de $I$ (par opposition \`a $I^r$, 
qui d\'esigne l'ensemble des \'el\'ements homog\`enes de degr\'e 
$r$ de $I$).

\subsubsection{D\'erivation des \'el\'ements de $S$}\label{derivation}
Soit $C_Y^*\subset \HH^0(Y,\anneau(1))^{\vee}$ le c\^one priv\'e de
l'origine de $Y\subset\PP\HH^0(Y,\anneau(1))^{\vee}$. Soit
$\pi:C_Y^*\to Y$ la projection naturelle. Pour tout $a\in\C^*$,
soit $h_a:C_Y^*\to C_Y^*$ l'homoth\'etie de rapport $a$. 
Soit ${\mathcal T}_Y$ le faisceau tangent de $Y$ et soit 
${\mathcal T}_{C_Y^*}$ le faisceau tangent de $C_Y^*$. Enfin, pour tout
$i\in\ZZ$ nous notons  
$\widetilde{T}^i\subset\HH^0(C_Y^*,{\mathcal T}_{C_Y^*})$
le sous-espace des sections $\widetilde{v}$ v\'erifiant
$(h_a)_*\widetilde{v}=a^{-i}\widetilde{v}$ pour tout $a\in\C^*$. 
Si nous identifions $S$ \`a un sous-anneau de l'anneau des fonctions
sur $C^*_Y$ alors le $S$-module gradu\'e $\widetilde{T}=\bigoplus_{i\in\ZZ}
\widetilde{T}^i$ d\'efinit un faisceau localement libre 
$\widetilde{\mathcal T}_Y$ sur $Y$ (si bien que 
$\HH^0(Y,\widetilde{\mathcal T}_Y\otimes\anneau_Y(i))=\widetilde{T}_i$)
qui s'ins\`ere dans la suite exacte
$$0\to\anneau_Y\to\widetilde{\mathcal T}_Y\to{\mathcal T}_Y\to 0.$$
La d\'eriv\'ee de Lie le long d'un champ de vecteurs sur $C_Y^*$ donne
pour tous les entiers relatifs $i$ et $j$ une application 
\begin{eqnarray*}
 \widetilde{T}^i\times S^j\to S^{i+j},
 \quad(\widetilde{v},P)\mapsto L_{\widetilde{v}}P.
\end{eqnarray*}

\subsubsection{\'Enonc\'e}
La proposition~\ref{3prop} ci-dessous 
est une adaptation de la proposition 3 de~\cite{NL} qui traite le cas
$Y=\PP^{N+1}$ et $n=N/2$. La premi\`ere assertion est une traduction 
alg\'ebrique de la transversalit\'e de Griffiths. La seconde est une 
cons\'equence de la formule de Leibniz pour la d\'erivation et s'inspire de
l'\'enonc\'e~(3.4.6) de~\cite{NLclaire2}. Enfin  la troisi\`eme traduit le
fait que les formes diff\'erentielles exactes ont une classe de
cohomologie nulle. La preuve de la proposition~\ref{3prop}
occupe la section~\ref{preuvegeo}.

\begin{prop} \label{3prop}
Pour tout $r\in\NN\cup\{-1\}$ nous avons, avec la convention $E_{-1}=S$,
\begin{enumerate}
\item $E_r=(E_{r+1}:F)$.
\item Pour un point g\'en\'erique $F$ de 
      $({\mathcal N}_{\lambda,U}^{d,n})^{\rm red}$ nous avons
      ${(E_r)}^2\tau_F
      \left(T_F\,({\mathcal N}_{\lambda,U}^{d,n})^{\rm red}\right)
      \subset E_{r+1}$.
\item Pour tout $P\in E_r^{\bullet}$ et pour
      tout $\widetilde{v}\in \widetilde{T}^{j}$ nous avons
      $$(L_{\widetilde{v}}P)F-(n+r+1)PL_{\widetilde{v}}F \in
      E_{r+1}^{\bullet+d+j}.$$
\end{enumerate}
\end{prop}

\subsubsection{Les id\'eaux $E_r$ et l'id\'eal jacobien}\label{jaco}
Le $S$-module $\widetilde{T}$ est de type fini, il est donc 
engendr\'e en degr\'e inf\'erieur \`a un entier $s\in\NN$ 
qui ne d\'epend que de $Y$ et de $\anneau_Y(1)$.
Soit $J_F$ l'id\'eal jacobien de $F$, engendr\'e par $F$ 
et par les d\'eriv\'ees partielles 
$L_{\widetilde{v}}F$, $\widetilde{v}\in \widetilde{T}$, donc
engendr\'e en degr\'e inf\'erieur  \`a $d+s$. Comme $X_F$ 
est lisse, l'id\'eal $J_F$ est sans points base. 
Les assertions~1 et~3 de la proposition~\ref{3prop} 
(pour $r=-1$ et P=1) impliquent donc

\begin{cor}
$J_F\subset E_0$~; en particulier $E_0^{d+s}$ est sans points-base.
\end{cor}

\subsection{\'Etude asympotique de la fonction de Hilbert des 
id\'eaux $E_r$} \label{alg}
Dans cette section nous consid\'erons le probl\`eme 
purement alg\'ebrique suivant. 

\subsubsection{Donn\'ee agl\'ebrique}\label{doal}
Soit $S$ l'alg\`ebre gradu\'ee d'un sch\'ema projectif lisse $Y$
de dimension $N+1$ muni d'un faisceau inversible tr\`es ample 
$\anneau_Y(1)$. Nous fixons des entiers
strictement positifs $u$, $s$, $m$, $b$, $d$ et $n\in\{1,\dots, N\}$. 
Nous consid\'erons les id\'eaux homog\`enes $I\subset S$  v\'erifiant 
\begin{enumerate}
\item $I^{d+s}$ est sans points base,
\item il existe un entier $D\geq (n+1)d-m$ et des id\'eaux Gorenstein
      $I_i$, $i\in\{1,\dots, u\}$ de degr\'e du socle $D$ tels que
      $I=\bigcap_{i=1}^{u}I_i$ et
\item  $\dim(S^d/I^d) \leq b\frac{d^n}{n!}$
\end{enumerate}

\subsubsection{Motivation}\label{moti}
Nous adoptons les notations de la section~\ref{defo}. Nous supposons   
de plus
$${\rm codim}\,({\mathcal N}_{\lambda,U}^{d,n},U) \leq b\frac{d^n}{n!}$$ 
comme dans l'hypoth\`ese du th\'eor\`eme~\ref{nlh+r}. L'\'enonc\'e du
th\'eor\`eme~\ref{nlh+r} \'etant stable par transport
plat, nous pouvons supposer sans perdre en g\'en\'eralit\'e
$F\in{\mathcal N}_{\lambda,U}^{d,n}$ g\'en\'erique et par cons\'equent
$$\lambda_F\not\in F^{n+1}\HH^N(X_F,\C)$$ 
(si $\lambda_F\in F^{n+1}\HH^N(X_F,\C)$ pour $F$ g\'en\'erique, nous 
aurions ${\rm codim}\,({\mathcal N}_{\lambda,U}^{d,n+1},U)\leq 
\frac{b(n+1)}{d}\frac{d^{n+1}}{(n+1)!}$, l'\'etude de 
${\mathcal N}_{\lambda,U}^{d,n+1}$ se heurterait 
\`a une contradiction pour $d>b(n+1)$ d\`es la proposition~\ref{alh}). 
Nous avons alors $E_r\neq S$ pour tout $r\in\NN$.

L'id\'eal $I=E_0$ v\'erifie alors la donn\'ee 
alg\'ebrique~\ref{doal} et les id\'eaux $E_r$ 
v\'erifient l'assertion~2 de la donn\'ee 
alg\'ebrique~\ref{doal}. 
En effet, l'assertion~1 pour $E_0$ est donn\'ee par le corollaire 
de la proposition~\ref{3prop}. L'assertion~2 pour $E_r$ r\'esulte 
par la section~\ref{defogo} de ce que $E_r\neq S$.
Enfin, l'assertion~3 pour $E_0$ r\'esulte des isomorphismes
$\tau_F:T_FU\simeq S^d$ et 
$\tau_F:T_F{\mathcal N}_{\lambda,U}^{d,n}\simeq E_0^d$.
Nous verrons aussi que les id\'eaux $E_r$ ne sont pas loin 
de v\'erifier l'assertion~3.

Les entiers $u$, $s$ et $m$ ainsi que 
l'alg\`ebre $S$ ne d\'ependent que de $Y$ et de $\anneau_Y(1)$. 

\subsubsection{\'Enonc\'e}\label{enoncealg}
L'assertion 1 de proposition~\ref{alh} ci-dessous g\'en\'eralise 
la proposition~1 de \cite{NL} qui traite le cas $Y=\PP^{N+1}_{\C}$ 
et $I$ Gorenstein. Elle repose sur la comparaison pour 
$k\to\infty$ entre la fonction de Hilbert $k\mapsto\dim\,S^k/I^k$ 
et la dimension des lieux base des syst\`emes lin\'eaires 
$I^k\subset S^k$.
L'assertion~2 repose sur une majoration asymptotique de la 
r\'egularit\'e des id\'eaux $\sqrt{I_V}$ et ${\mathfrak q}$
au sens de Bayer et Mumford. La preuve de la proposition~\ref{alh}
occupe la section~\ref{preuvealg}.

\begin{prop}\label{alh}
Pour tout $\varepsilon \in \RR^*_+$ il existe une constante
$C\in\RR_+^*$ qui ne d\'epend que de $\varepsilon$, 
$S$, $u$, $s$ et $m$, telle que pour tout $d \geq Cb$ et pour tout 
id\'eal homog\`ene $I\subset S$  v\'erifiant la donn\'ee 
alg\'ebrique~\ref{doal} l'assertion suivante est vraie.

Il existe un sch\'ema $V\subset Y$ de dimension pure $n$ et de
degr\'e inf\'erieur  \`a $(1+\varepsilon)b$ tel que
\begin{enumerate}
\item l'id\'eal $I_V$ de $V$ est inclus dans $I$
\item les id\'eaux $I_W$ et ${\mathfrak q}$ sont engendr\'es en 
degr\'e inf\'erieur  \`a $Cb^a$, 
o\`u $I_W$ d\'esigne l'id\'eal de $W=V^{\rm red}$, ${\mathfrak q}$ 
est d\'efini ci-dessous et $a=2^{\hh^0(Y,\anneau(1))-1}$.  
\end{enumerate}
\end{prop}

\subsubsection{D\'efinition de l'id\'eal ${\mathfrak q}$}\label{defq}
Comme $V$ est de dimension pure $n$, $I_V$ est une
intersection finie d'id\'eaux primaires de hauteur $n$ not\'es 
${\mathfrak P}_i$, et $I_W$ est une intersection finie d'id\'eaux
premiers not\'es ${\mathfrak p}_i$, o\`u 
${\mathfrak p}_i=\sqrt{{\mathfrak P}_i}$.
Pour tout $i$, nous notons $b_i$ la longueur de l'anneau artinien
$(S/I_V)_{{\mathfrak p}_i}= (S/{\mathfrak P}_i)_{{\mathfrak p}_i}$ et 
$\delta_i$ le degr\'e du sch\'ema $V({\mathfrak p}_i)$. Nous avons
$({\mathfrak p}_i)^{b_i}\subset {\mathfrak P}_i$~;  
soit $\beta_i< b_i$ le plus grand entier tel que 
$({\mathfrak p}_i)^{\beta_i}\not\subset {\mathfrak P}_i$. 
Nous posons ${\mathfrak q}=\prod_i({\mathfrak p}_i)^{\beta_i}$ et
$\beta=\sum_i\beta_i\delta_i$. Nous avons 
$\beta\leq\sum_i\delta_i(b_i-1)=
\deg V-\deg V^{\rm red}\leq (1-\varepsilon)b$.

Nous avons ${\mathfrak q}I_W\subset I_V$ et 
si $V$ est r\'eduit, alors ${\mathfrak q}=S$.

\subsection{Les id\'eaux $E_r$ et le sch\'ema $V$}\label{algeco}
Les notations sont celles de la section~\ref{defo}. 
Nous supposons satisfaites les hypoth\`eses de la section~\ref{moti}.
Comme $F\in{\mathcal N}_{\lambda,U}^{d,n}$ est g\'en\'erique
nous pouvons supposer sans perdre en g\'en\'eralit\'e que
$F$ v\'erifie l'assertion~2 de la proposition~\ref{3prop}.

L'id\'eal $E_0$ v\'erifie la donn\'ee alg\'ebrique~\ref{doal}.
Nous fixons $\varepsilon\in\RR_+^*$ et nous supposons $d\geq Cb^a$, 
o\`u $a=2^{\hh^0(Y,\anneau(1))-1}$
et o\`u $C$ est la constante donn\'ee par la proposition~\ref{alh}. 
Nous notons $V$ et $W=V^{\rm red}$ les sch\'emas associ\'es \`a $E_0$ 
par la proposition~\ref{alh}.

La proposition~\ref{I^2am} et la proposition~\ref{Zred}
ci-dessous g\'en\'eralisent alors respectivement le lemme~2 et le
th\'eor\`eme~2 de~\cite{NL}. Les preuves occupent respectivement 
les sections~\ref{preuveI^2am} et~\ref{preuveZred}.

\begin{prop} \label{I^2am}
Il existe une constante $C\in\RR_+^*$ qui ne d\'epend que de 
$\varepsilon$, de $Y$ et de $\anneau_Y(1)$ telle que pour tout  
$d\geq Cb^a$ et pour tout $r\in\NN$ nous avons 
${\mathfrak q}(I_W)^{r+1}\subset {E_{r}}$.
\end{prop}

\begin{prop} \label{Zred} Supposons $n\leq [N/2]$. 
Il existe une constante $C\in\RR_+^*$ qui ne d\'epend que de 
$\varepsilon$, de $Y$ et de $\anneau_Y(1)$ telle que pour tout  
$d\geq Cb^a$ nous avons $F\in I_W$.
\end{prop}

\section{Description globale des vari\'et\'es de Hodge}\label{glo}

\subsection{Notations et rappels de topologie}\label{rappelstopo}
Les notation sont celle de la section~\ref{introtopo}.

Le groupe $\pi_1(\YYY,F)$ agit par monodromie sur $\HH^{N}(X_F,\C)$. 
Plus pr\'ecis\'ement, la th\'eorie de Lefschetz affirme que la
repr\'esentation de monodromie de $\pi_1(\YYY,F)$ sur $\HH^{N}(X_F,\C)$
est somme directe de la repr\'esentation triviale sur 
$j_F^*\HH^{N}(Y,\C)$ et d'une repr\'esentation irr\'eductible sur
$\HH^{N}(X_F,\C)_{\ev}$.  

Si $W$ est une sous-vari\'et\'e de $X_F$ de dimension $n$, 
nous notons $\YYY(W)\subset\YYY$ l'espace param\'etrant les
hypersurfaces lisses contenant $W$. Nous \'etudions l'action de
monodromie de $\pi_1(\YYY(W),F)$ sur $\HH^{N}(X_F,\C)$~: 
pour tout sous-ensemble $E\subset\HH^{N}(X_F,\C)_{\ev}$ nous notons 
$H(E)\subset\HH^{N}(X_F,\C)_{\ev}$ la $\pi_1(\YYY(W),F)$-repr\'esentation 
engendr\'ee par $E$ et $E^{\perp}\subset\HH^{N}(X_F,\C)_{\ev}$ 
l'orthogonal de $E$ pour le cup-produit. 

\subsection{\'Etude de $H(\lambda_F)$}\label{res}
Nous adoptons les notations et les hypoth\`eses des 
sections~\ref{algeco} et~\ref{rappelstopo}. 
Nous supposons $d\geq Cb^a$, o\`u $C\in\RR^*_+$ est une 
constante v\'erifiant les propositions~\ref{alh},~\ref{I^2am} 
et~\ref{Zred}. 

L'id\'ee de la proposition~\ref{resume} ci-dessous est de consid\'erer 
les \'equations de ${\mathcal N}_{\lambda,U}^{d,n}$ donn\'ees 
par les \'el\'ements de degr\'e $nd$ de ${\mathfrak q}K$ 
(voir la section~\ref{equa}). 
La proposition~\ref{I^2am} implique alors que les 
d\'eriv\'ees en $F$ \`a tout ordre le long de $T_F\YYY(W)$
de ces \'equations sont nulles. Par unicit\'e du prolongement 
analytique nous en d\'eduisons une condition sur les images 
par transport plat dans $\YYY(W)$ de la classe $\lambda_F$.
 
\begin{prop}\label{resume} 
$H(\lambda_F)\subset\phi_{F,n}({\mathfrak q}K^{nd})^{\perp}$.
\end{prop}

\preuve
Soient $Q\in{\mathfrak q}$ et $\omega\in K$ des \'el\'ements homog\`enes 
dont la somme des degr\'es est ${nd}$.
Comme $\lambda_F\in F^n\HH^{N}(X_F,\C)$, nous avons 
$\lambda_F\smile\phi_{F,n}(Q\omega)=0$. 

\begin{lemme}\label{resumelocal}
Soit $U_W$ un voisinage simplement connexe de $F$ dans $\YYY(W)$. 
Pour tout $G\in U_W$ nous avons 
$\lambda_G\smile\phi_{G,n}(Q\omega)=0$, o\`u $\lambda_G$ d\'esigne 
l'image de $\lambda_F$ dans $\HH^N(X_{G},\C)$ par transport 
plat dans $U_W$.
\end{lemme}

\preuve
La fonction $U_W\to\C$, $G\mapsto\lambda_G\smile\phi_{G,n}(Q\omega)$ 
est analytique et s'annule en $G=F$.
Par unicit\'e du prolongement analytique, il suffit de montrer
que pour tout $r\in\NN$ et pour tout $\vec{v}\in T_FU_W$ nous avons
$\left(\frac{\partial}{\partial\vec{v}}\right)^{r+1}
(G\mapsto\lambda_G\smile\phi_{G,n}(Q\omega))(F)=0.$

Or d'apr\`es la section~\ref{alte}, nous avons
\begin{eqnarray*}
   \left(\frac{\partial}{\partial\vec{v}}\right)^{r+1}
   (G\mapsto\lambda_G\smile\phi_{G,n}(Q\omega))(F)
&=&\left(\frac{\partial}{\partial\vec{v}}\right)^{r+1}
   \left(G\mapsto\sum_{i=1}^t\lambda_i\int_{{\rm Tub}\,\gamma_i} 
   \frac{Q\omega}{G^n}\right)(F)\\
&=&\frac{(-1)^{r+1}(n+r)!}{(n-1)!}
   \sum_{i=1}^t\lambda_i\int_{{\rm Tub}\,\gamma_i} 
   \frac{\tau_F(\vec{v})^{r+1}Q\omega}{F^{r+n+1}}.
\end{eqnarray*}
Comme $\vec{v}\in T_FU_W$, donc $\tau_F(\vec{v})\in I_W$ et donc 
$(\tau_F(\vec{v}))^{r+1}Q\in(I_W)^{r+1}{\mathfrak q}$, 
d'apr\`es la proposition~\ref{I^2am} nous avons 
$(\tau_F(\vec{v}))^{r+1}Q\in E_{r}$, donc 
$\left(\frac{\partial}{\partial\vec{v}}\right)^{r+1}
(G\mapsto\lambda_G\smile\phi_{G,n}(Q\omega))(F)=0$.
\cucu

\bigskip

\noindent{\it Fin de la preuve de la proposition~\ref{resume}.\ ---\ }
Soit $p_W:\widetilde{\YYY(W)}\to\YYY(W)$ le rev\^etement universel de
$\YYY(W)$ et soit $\widetilde{F}\in p_W^{-1}(F)$. Pour tout
$G\in\YYY(W)$ et $\widetilde{G}\in p_W^{-1}(G)$ nous notons
$\lambda^{\widetilde{F},\widetilde{G}}\in\HH^N(X_G,\C)$ l'image de 
$\lambda_F$ par transport plat le long de l'image par $p_W$ 
d'un chemin reliant $\widetilde{F}$ \`a $\widetilde{G}$. Le lieu
$\widetilde{{\mathcal N}_{\lambda,\mu}}=
\{\widetilde{G}\in\widetilde{\YYY(W)}\ |\ 
\lambda^{\widetilde{F},\widetilde{G}}\smile\phi_{G,n}(Q\omega)=0\}$ 
est un sous-ferm\'e analytique de $\widetilde{\YYY(W)}$~; 
comme $\widetilde{\YYY(W)}$ est
irr\'eductible, le lemme~\ref{resumelocal} implique 
$\widetilde{{\mathcal N}_{\lambda,\mu}}=\widetilde{\YYY(W)}$.
En particulier, $\widetilde{{\mathcal N}_{\lambda,\mu}}$ contient
$p_W^{-1}(F)$, ce qui implique 
$H(\lambda_F)\subset\phi_{F,n}(Q\omega)^{\perp}$.
\cucu

\medskip

\noindent
{\it Remarque.\ ---\ }\
Lorsque $V$ est r\'eduit la proposition~\ref{resume} donne 
$H(\lambda_F)\subset F^n\HH^N(X_F,\C)_{\ev}$. 
Plus pr\'ecis\'ement, la preuve implique qu'alors $\YYY(W)$ est inclus 
dans ${\mathcal N}^{d,n}_{\lambda,U}$. 

\subsection{Un \'enonc\'e de monodromie}\label{monopre}
Nous adoptons les notations des sections~\ref{rappelstopo} 
et~\ref{introtopo}.
Pour toute sous-vari\'et\'e $W$ de $X_F$ de dimension $n$
l'espace $H(W)$ est $\pi_1(\YYY(W),F)$-invariant, et
comme $H(W)\subset\HH^{n,n}(X_F)_{\ev}$, 
le th\'eor\`eme d'indice de Hodge implique que la restriction du 
cup-produit \`a $H(W)$ est non d\'eg\'en\'er\'ee et nous avons
l'\'egalit\'e de $\pi_1(\YYY(W),F)$-modules 
$\HH^{N}(X_F,\C)_{\ev}=H(W)\oplus H(W)^{\perp}$.

Nous supposons $n\in\{1,\dots,[N/2]\}$ et nous notons $e\in\NN$ 
le plus petit entier tel que $W$ est engendr\'e en degr\'e 
strictement inf\'erieur \`a $e$. Remarquons que lorsque $W$ est 
donn\'e par la proposition~\ref{alh}, nous avons $e\leq Cb^a+1$.   
Nous avons montr\'e dans \cite{topo}~: 

\begin{prop}\label{mono}
Il existe une constante $C\in\RR^*_+$ qui ne d\'epend que de 
$Y$ et de $\anneau_Y(1)$, telle que pour tout $d\geq Ce$
la $\pi_1(\YYY(W),F)$-repr\'esentation de monodromie 
sur $\HH^{N}(X_F,\C)_{\ev}$ est
\begin{enumerate}
\item irr\'eductible pour $N>2n$~;
\item somme directe d'une repr\'esentation triviale sur $H(W)$ et d'une
      repr\'esentation irr\'eductible sur $H(W)^{\perp}$ pour $N=2n$.
\end{enumerate}
\end{prop}

\subsection{Preuve du th\'eor\`eme~\ref{nlh+r}}
Nous adoptons les hypoth\`eses du th\'eor\`eme~\ref{nlh+r}.
Nous fixons $\varepsilon\in\RR_+^*$ et
nous supposons $d\geq Cb^a$, o\`u $C\in\RR^*_+$ est une 
constante v\'erifiant les propositions~\ref{resume} et~\ref{mono}. 

Nous avons $\phi_{F,n}({\mathfrak q}K^{nd})\neq 0$
(en effet, si nous choisissons par exemple
$\omega\in K^{d-\beta}$ et  
$Q\in{\mathfrak q}^{\beta}$ tels que la section 
$Q\omega \in K^d$ ne s'annule pas sur $X_F$, la classe
$\phi_{F,n}(F^{n-1}Q\omega)=
{\rm res}\left[\frac{Q\omega}{F}\right]$ est non nulle),
donc $\phi_{F,n}({\mathfrak q}K^{nd})^{\perp}\neq\HH^{N}(X_F,\C)_{\ev}$.

Si $2n<N$, l'assertion~1 de la proposition~\ref{mono} implique
$H(\lambda_F)=\HH^N(X_F,\C)_{\ev}$, or la proposition~\ref{resume}
implique $\phi_{F,n}({\mathfrak q}K^{nd})^{\perp}=\HH^N(X_F,\C)_{\ev}$, 
contradiction~: la classe $\lambda_F$ n'existe pas.

Si $2n=N$, supposons par l'absurde $\lambda_F\not\in H(W)$. Il
existe alors des classes $\gamma_F\in H(W)$ et
$\nu_F\in H(W)^{\perp}\setminus\{0\}$ telles que 
$\lambda_F=\gamma_F+\nu_F$. D'une part, comme $\dim\,H(W)^{\perp}>1$,  
l'assertion~2 de la proposition~\ref{mono} implique 
que la $\pi_1(\YYY(W),F)$-repr\'esentation sur $H(W)^{\perp}$ est
irr\'eductible et non triviale, donc
$H(\lambda_F)=\C\gamma_F\oplus H(W)^{\perp}$, ce qui implique 
$H(W)^{\perp}\subset 
H(\lambda_F)\subset \phi_{F,n}({\mathfrak q}K^{nd})^{\perp}$ 
d'apr\`es la proposition~\ref{resume}.
D'autre part nous avons 
$H(W)\subset F^n\HH^N(X_F,\C)_{\ev}\subset 
\phi_{F,n}({\mathfrak q}K^{nd})^{\perp}.$ Nous avons donc 
$H(W)^{\perp}\oplus H(W)\subset \phi_{F,n}({\mathfrak q}K^{nd})^{\perp}$,
donc $\phi_{F,n}({\mathfrak q}K^{nd})^{\perp}=\HH^N(X_F,\C)_{\ev}$, 
contradiction. Donc $\lambda_F\in H(W)$. 
\cucu
\section{Preuve de la proposition~\ref{3prop}}\label{preuvegeo}

L'assertion~1 de la proposition~\ref{3prop} est \'evidente. 

\subsection{Preuve de l'assertion~2}
Pour tout $G\in{\mathcal N}_{\lambda,U}^{d,n}$ d\'efinissons 
l'id\'eal $E_r(G)$ associ\'e \`a la classe $\lambda _G$, par analogie
avec l'id\'eal $E_r=E_r(F)$ associ\'e \`a la classe $\lambda _F$.
Pour tout ${k}\in\NN$ posons
$${\mathcal E}_r^{k}=\left\{(P,G)\in S^{k}\times 
({\mathcal N}_{\lambda,U}^{d,n})^{\rm red}\ |\ P\in E_r^{k}(G)\right\}.$$
Nous choisissons un point lisse 
$F\in({\mathcal N}_{\lambda,U}^{d,n})^{\rm red}$ 
tel que pour tout $k\in\NN$ il existe un voisinage lisse 
$U_{r,k}$ de $F$ dans $({\mathcal N}_{\lambda,U}^{d,n})^{\rm red}$ tel que
la restriction de la projection
$pr_{k}:{\mathcal E}_r^{k}\to({\mathcal N}_{\lambda,U}^{d,n})^{\rm red}$ 
\`a $pr_{k}^{-1}(U_{r,k})$ est un fibr\'e trivial. 

Soient $i$ et $j$ deux entiers positifs, $P\in E_r^i$, $Q\in E_r^j$
et $R=\tau_F(\vec{r})$, avec 
$\vec{r}\in T_F\,({\mathcal N}_{\lambda,U}^{d,n})^{\rm red}$~: nous
devons montrer $PQR\in E_{r+1}$. 

Posons $U_{r,i,j}=U_{r,i}\cap U_{r,j}$~: nous avons des fibr\'es
triviaux $pr_{i}:pr_{i}^{-1}(U_{r,i,j})\to U_{r,i,j}$ et
$pr_{j}:pr_{j}^{-1}(U_{r,i,j})\to U_{r,i,j}$, ils poss\`edent donc
des sections ${\bf P}:U_{r,i,j}\to pr_{i}^{-1}(U_{r,i,j})$ et 
${\bf Q}:U_{r,i,j}\to pr_{j}^{-1}(U_{r,i,j})$, telles que
${\bf P}(F)=P$ et ${\bf Q}(F)=Q$~; nous identifions ces sections \`a
des fonctions analytiques    
${\bf P}:U_{r,i,j}\to S^i$ et ${\bf Q}:U_{r,i,j}\to S^j$
telles que ${\bf P}(G)\in E_r^i(G)$ et
${\bf Q}_i(G)\in E_r^j(G)$ pour tout $G\in U_{r,i,j}$. 

Pour tout   
$\omega\in K^{(n+r+1)d-i-j}$ la fonction
$$U_{r,i,j}\to\C,\qquad G\mapsto 
  \sum_{i=1}^t\lambda_i\int_{{\rm Tub}\,\gamma_i}
  \frac{{\bf P}(G){\bf Q}(G)\omega}{{G}^{n+r+1}}$$
est identiquement nulle~; en la d\'erivant en $F$ le long du vecteur
$\vec{r}$ nous obtenons l'\'egalit\'e
$$\sum_{i=1}^t\lambda_i\int_{{\rm Tub}\,\gamma_i}\left(
  \frac{\frac{\partial{\bf P}}{\partial\vec{r}}(F)Q\omega}{{F}^{n+r+1}}
+ \frac{P\frac{\partial{\bf Q}}{\partial\vec{r}}(F)\omega}{{F}^{n+r+1}}
- (n+r+1)\frac{RPQ\omega}{{F}^{n+r+2}}\right)=0.$$
Comme $P\in E_r^i$ et $Q\in E_r^j$, les deux premiers termes sont
nuls, ce qui implique l'assertion~2.
\cucu

\subsection{Preuve de l'assertion~3}
Notons ${\mathcal K}_{C_Y^*}$ le faisceau canonique de $C_Y^*$. Pour
tout $i\in\ZZ$ notons 
$\widetilde{K}^i\subset\HH^0(C_Y^*,{\mathcal K}_{C_Y^*})$
le sous-espace des sections $\omega$ v\'erifiant $h_a^*\omega=a^i\omega$
pour tout $a\in\C^*$. Si $\vec{e}\in\widetilde{T}^1$ d\'esigne le
champ d'Euler (radial et invariant par homoth\'etie), nous avons pour
tout $i\in\ZZ$ des isomorphismes canoniques 
${\rm Int}_{\vec{e}}:\widetilde{K}^i\to K^i.$

Soient $P\in E_r$, $\widetilde{v}\in\widetilde{T}$ et
$\omega\in\bigoplus_{j\in\ZZ}\widetilde{K}_j$
des \'el\'ements homog\`enes dont la somme des degr\'es est \'egale
\`a  $(n+r+1)d$.
Nous avons l'\'egalit\'e des $(N+1)$-formes diff\'erentielles 
homog\`enes de degr\'e $0$ sur $C_Y^*\setminus \pi^{-1}(X_F)$~: 
\begin{eqnarray*}\label{cartan}
L_{\widetilde{v}}\left(\frac{P\omega}{F^{n+r+1}}\right)=
\frac{PL_{\widetilde{v}}\,\omega}{F^{n+r+1}}
+\frac{(L_{\widetilde{v}}P)\omega}{F^{n+r+1}}
-(n+r+1)\frac{(L_{\widetilde{v}}F)P\omega}{F^{n+r+2}}.
\end{eqnarray*}
Nous \'etudions l'image par ${\rm Int}_{\vec{e}}$ de cette
\'egalit\'e, que nous interpr\'etons comme une \'egalit\'e entre formes
diff\'erentielles sur $Y\setminus X_F$.

En appliquant deux fois la formule de Cartan-Lie 
$L=D\circ{\rm Int}+{\rm Int}\circ D$ \`a l'image par
${\rm Int}_{\vec{e}}$ du membre de gauche, nous obtenons
successivement~:    
\begin{eqnarray*}\label{eqcartan}
{\rm Int}_{\vec{e}}\,L_{\widetilde{v}}
\left(\frac{P\omega}{F^{n+r+1}}\right)
={\rm Int}_{\vec{e}}\,D\,{\rm Int}_{\widetilde{v}}
\left(\frac{P\omega}{F^{n+r+1}}\right)
=-D\,{\rm Int}_{\vec{e}}\,{\rm Int}_{\widetilde{v}}
\left(\frac{P\omega}{F^{n+r+1}}\right)
\end{eqnarray*}
(la premi\`ere \'egalit\'e r\'esulte de ce que 
$D\left(\frac{P\omega}{F^{n+r+1}}\right)=0$, puisque $\omega$ est de
degr\'e maximal, et la seconde de ce que $L_{\vec{e}}\,
{\rm Int}_{\widetilde{v}}\left(\frac{P\omega}{F^{n+r+1}}\right)=0$, puisque
la forme diff\'erentielle ${\rm Int}_{\widetilde{v}}
\left(\frac{P\omega}{F^{n+r+1}}\right)$ est invariante par
homoth\'etie). En particulier, la forme 
${\rm Int}_{\vec{e}}\,L_{\widetilde{v}}
\left(\frac{P\omega}{F^{n+r+1}}\right)$ est exacte.

Nous avons donc 
$$\sum_{i=1}^t\lambda_i\int_{{\rm Tub}\,\gamma_i}\left(
\frac{P{\rm Int}_{\vec{e}}\,L_{\widetilde{v}}\,\omega}{F^{n+r+1}}
+\frac{(L_{\widetilde{v}}P){\rm Int}_{\vec{e}}\,\omega}{F^{n+r+1}}
-(n+r+1)\frac{(L_{\widetilde{v}}F)P{\rm Int}_{\vec{e}}\,\omega}{F^{n+r+2}}
\right)=0.$$
Comme $P\in E_r^{\bullet}$, le premier terme est nul, donc
$$\sum_{i=1}^t\lambda_i\int_{{\rm Tub}\,\gamma_i}\frac
{(L_{\widetilde{v}}P)F{\rm Int}_{\vec{e}}\omega-
(n+r+1)(L_{\widetilde{v}}F)P{\rm Int}_{\vec{e}}\omega}
{F^{n+r+2}}=0.$$ 
Ceci implique l'assertion~3, par d\'efinition de l'id\'eal $E_{r+1}$.
\cucu

\section{Preuve de la proposition~\ref{alh}}\label{preuvealg} 

\subsection{Notations et rappels}\label{rappelsalg}
Les notations sont celles de la section~\ref{doal}. Nous posons
$A=\dim\,S^1=\hh^0(Y,\anneau(1))$ et nous notons $\mathfrak S$ 
l'alg\`ebre des polyn\^omes \`a $A$ variables. Le plongement
$Y\subset \PP^{A-1}_{\C}$ donn\'e par 
$\anneau_Y(1)$ identifie $S$ \`a un quotient de $\mathfrak S$.   

\subsubsection{Lieu-base d'un id\'eal}    
Soit $I\subset S$ un id\'eal homog\`ene.
Pour tout $i\in\{0,\dots,N+1\}$ nous notons $l_i(I)\in \NN \cup \infty$
le plus petit nombre $l\in \NN \cup \infty$ tel que le lieu-base de
$I^{l}$ est de dimension inf\'erieure ou \'egale \`a 
$N-i$ (en convenant que la partie vide est de dimension $-1$). De
mani\`ere \'equivalente $l_i(I)$ est le plus grand nombre
$l\in\NN\cup\infty$ tel que le lieu-base de $I^{l-1}$ est de dimension
sup\'erieure ou \'egale \`a $N-i+1$. Nous avons $l_0(I)\leq \cdots \leq
l_{N+1}(I)$. L'id\'eal $I$ est sans points base si et seulement si
$l_{N+1}(I)<\infty$~; le cas \'ech\'eant l'id\'eal $I$ contient une
intersection compl\`ete de multidegr\'e $(l_0(I),\dots,l_{N+1}(I))$. 

\subsubsection{Id\'eaux Gorenstein}\label{gore}
Nous disons que l'id\'eal $I\subset S$ est Gorenstein de degr\'e du
socle $D$ si $I^{D+1}=S^{D+1}$ et s'il existe une forme lin\'eaire 
$\Lambda\in(S^D)^\vee$ telle que pour tout $d\in\{0,\dots,D\}$ nous
avons 
$$I^d =\{P\in S^d\ |\ \forall Q\in S^{D-d},\ \Lambda(PQ)=0\}.$$
Remarquons que l'id\'eal $I\subset S$ est Gorenstein si et seulement
si sa pr\'eimage dans ${\mathfrak S}$ l'est. La forme lin\'eaire
$\Lambda$ induit alors pour tout $d \in \{ 0, \dots ,D\}$ 
un isomorphisme canonique $S^d/I^d \simeq (S^{D-d}/I^{D-d})^{\vee}$. 

Soient $I$ et $I'$ deux id\'eaux Gorenstein associ\'es aux formes
lin\'eaires $\Lambda \in (S^D)^\vee$ et $\Lambda' \in (S^{D'})^\vee$.
Supposons $I \subset I'$ (donc $D \geq D'$). 
Alors l'image de $\Lambda'$ dans $(S^{D'}/I^{D'})^\vee$ s'identifie via
l'isomorphisme induit par $\Lambda$ \`a un polyn\^ome $\bar{F} \in
S^{D-D'}/I^{D-D'}$. Pour tout polyn\^ome $Q \in S^{D'}$, nous avons 
$\Lambda '(Q) = \Lambda (QF)$, o\`u $F\in S^{D-D'}$ est un relev\'e de
$\bar{F}$. Nous avons donc $I'=(I:F)$. 

\subsubsection{Intersections compl\`etes}\label{subalge}
Soit $I\subset S$ un id\'eal homog\`ene intersection compl\`ete 
sans points base. Lorsque $S$ est une alg\`ebre de polyn\^omes, 
l'id\'eal $I$ est Gorenstein de degr\'e du socle 
$\sum_{i=0}^{N+1}(l_i(I)-1)$. Pour $S$ quelconque l'id\'eal $I$
n'est pas en g\'en\'eral Gorenstein. Cependant, il existe toujours
une constante $C$ qui ne d\'epend que de $S$ telle que $I$ 
co\"{\i}ncide avec $S$ en degr\'e $\sum_{i=0}^{N+1}l_i(I)+C$.

\subsection{G\'en\'eralisation des r\'esultats de~\cite{NL}}
Les notations sont celles des sections~\ref{doal} et~\ref{rappelsalg}.
 
Les lemmes~\ref{llh} et~\ref{alh-} ci-dessous g\'en\'eralisent
respectivement le lemme 1 et la proposition 1 de~\cite{NL} auquels
leurs preuves se ram\`ement.
Ces lemmes serviront dans les sections~\ref{preuveI^2am}
et~\ref{preuveZred}. L'assertion~1 du lemme~\ref{alh-} est 
l'assertion~1 de la proposition~\ref{alh}. Les assertions~2 et~3 du
lemme~\ref{alh-} serviront \`a la preuve de l'assertion~2 de la 
proposition~\ref{alh}, qui occupe la section~\ref{sebamu}. 

\begin{lemme}\label{llh}
Pour tout $\varepsilon \in \RR^*_+$ il existe une constante
$C\in\NN^*$ qui ne d\'epend que de $\varepsilon$, 
$S$, $u$ et $m$, telle que pour tout $d \geq Cb$ 
et pour tout id\'eal homog\`ene $I\subset S$  v\'erifiant les 
assertions~2 et~3 de la donn\'ee alg\'ebrique~\ref{doal},
nous avons $l_j(I)\leq\varepsilon d$ pour tout $j\in\{0,\dots,N-n\}$. 
\end{lemme}

\preuve
Dans le cas o\`u l'id\'eal $I$ est Gorenstein, la preuve du 
lemme~\ref{llh} est calqu\'ee sur celle du lemme 1 de~\cite{NL}. Le cas
g\'en\'eral s'en d\'eduit~: consid\'erons l'id\'eal $\prod_{i=1}^u I_i\subset
I$~: pour tout $j\in\{0,\dots,N-n\}$ nous avons
$l_j(I)\leq l_j(\prod_{i=1}^u I_i)\leq\sum_{i=1}^u l_j(I_i)$~; le
lemme~\ref{llh} pour les id\'eaux Gorenstein $I_i$, $i\in\{1,\dots,
u\}$ et pour $\frac{\varepsilon}{u}$ au lieu de $\varepsilon$ donne
$l_j(I_i)\leq\frac{\varepsilon}{u}d$, donc $l_j(I)\leq\varepsilon d$.  
\cucu

\begin{lemme} \label{alh-}  
Pour tout $\varepsilon \in \RR^*_+$ il existe une constante
$C\in\NN^*$ qui ne d\'epend que de $\varepsilon$, 
$S$, $u$, $s$ et $m$, telle que pour tout $d \geq Cb$
et pour tout id\'eal homog\`ene $I\subset S$  v\'erifiant 
la donn\'ee alg\'ebrique~\ref{doal}, il existe un 
sch\'ema $V\subset Y$ de dimension pure $n$ et de
degr\'e inf\'erieur  \`a $(1+\varepsilon)b$,
d'id\'eal $I_V\subset S$, tel que les assertions suivantes sont vraies
\begin{enumerate}
\item $I_V\subset I$,
\item ${\mathfrak I}_V=({\mathcal F}_1,\dots,{\mathcal F}_{A-n-1}:{\mathcal G})$, o\`u 
${\mathfrak I}_V$ d\'esigne la pr\'eimage de $I_V$ dans
 ${\mathfrak S}$, o\`u les ${\mathcal F}_i$, $i\in\{1,\dots,A-n-1\}$ sont des
polyn\^omes homog\`enes de ${\mathfrak S}$ de degr\'e au plus 
$\deg V$ en intersection compl\`ete, et  
o\`u ${\mathcal G}$ est un polyn\^ome homog\`ene de ${\mathfrak S}$ 
de degr\'e au plus $(A-n-1)\deg V+C$, 
\item les id\'eaux $I$ et $I_V$ co\"{\i}ncident en degr\'e 
inf\'erieur  \`a 
$d-(A-n-1)\deg V-C$.  
\end{enumerate}
\end{lemme}

\preuve
Dans le cas o\`u l'id\'eal $I$ est Gorenstein, la preuve du 
lemme~\ref{alh-} est calqu\'ee sur celle de la proposition~1
de~\cite{NL}. Nous allons en d\'eduire le cas g\'en\'eral.

Pour tout $i\in\{1,\dots,u\}$ 
l'id\'eal $I_i$ est d\'efini par une forme lin\'eaire
$\Lambda_i\in (S^{D})^{\vee}$. Soit $L\subset (S^{D})^{\vee}$ 
l'espace vectoriel engendr\'e par les $\Lambda_i$~; posons 
$L^*=L\setminus\{0\}$. 
Pour tout $\Lambda\in L^*$, l'id\'eal $I_\Lambda$ associ\'e \`a $\Lambda$ 
est Gorenstein et v\'erifie la donn\'ee alg\'ebrique~\ref{doal}~; 
soit $V_\Lambda$ le sch\'ema associ\'e.
Posons $V=\bigcup_{i=1}^{u}V_{\Lambda_i}$. C'est un sch\'ema de 
dimension pure $n$ et comme $I=\bigcap_{i=1}^{u}I_i$, son id\'eal
co\"{\i}ncide avec $I$ en degr\'e inf\'erieur  \`a 
$d-(A-n-1)\deg\,V-C$.

Comme $I=\bigcap_{\Lambda\in L^*}I_\Lambda$, nous avons 
$V=\bigcup_{\Lambda\in L^*}V_\Lambda$. 
Pour tout sous-sch\'ema strict $V'$ de $V$, l'espace 
$\{\Lambda\in L^*\ |\ V_\Lambda\subset V'\}\cup\{0\}$ est
donc un sous-espace vectoriel strict de $L$~; pour $\Lambda\in L^*$
n'appartenant pas \`a un nombre fini de sous-espaces vectoriels stricts
param\'etr\'es par les sous-sch\'emas stricts de $V$ de dimension pure
$n$ nous avons donc $V_\Lambda=V$. Le lemme~\ref{alh-} pour
$I_{\Lambda}$ implique alors le lemme~\ref{alh-} pour $I$. 
\cucu

\subsection{Les id\'eaux $I_W$ et ${\mathfrak q}$ sont
               engendr\'es en petit degr\'e}\label{sebamu} 
Dans cette section nous montrons l'assertion~2 de la 
proposition~\ref{alh}. Nous supposons v\'erifi\'ees les hypoth\`eses
num\'eriques de la proposition~\ref{alh}. Les notations sont celles de 
la section~\ref{defq}. 

\begin{lemme}\label{reduction}
\begin{enumerate}
\item Pour tout $j\in\{0,\dots,N-n\}$ nous avons 
$l_j({\mathfrak q})\leq\beta$.
\item Pour $Q\in {\mathfrak q}^{\beta}$ 
g\'en\'erique nous avons $I_W=(I_V:Q)$. 
\end{enumerate} 
\end{lemme}

\preuve
Comme $V({\mathfrak p}_i)$ est r\'eduit et que l'id\'eal
${{\mathfrak p}_i}$ contient en degr\'e $\delta_i$ les 
\'equations des c\^ones de base $V({\mathfrak p}_i)$, 
les sch\'emas $V({\mathfrak p}_i)$ et  
$V\left( {{\mathfrak p}_i}^{\delta_i}\right)$ 
ont m\^eme point g\'en\'erique. 

L'espace $\prod_i\left({{\mathfrak p}_i}^{\delta_i}\right)^{\beta_i}$ 
a donc pour lieu-base le sch\'ema 
$\bigcup_i V\left( {{\mathfrak p}_i}^{\delta_i}\right)$
de dimension $n$. Or
$\prod_i\left({{\mathfrak p}_i}^{\delta_i}\right)^{\beta_i}
\subset{\mathfrak q}^{\beta}$, donc
$V({\mathfrak q}^{\beta})\subset
\bigcup_i V\left( {{\mathfrak p}_i}^{\delta_i}\right)$
ce qui implique l'assertion~1.

Comme ${{\mathfrak p}_i}S_{{\mathfrak p}_i}=
{{\mathfrak p}_i}^{\delta_i}S_{{\mathfrak p}_i}$ et 
$({\mathfrak p}_i)^{\beta_i}S_{{\mathfrak p}_i}
\not\subset{{\mathfrak P}_i}S_{{\mathfrak p}_i}$, nous avons
${\left({{\mathfrak p}_i}^{\delta_i}\right)}^{\beta_i}
S_{{\mathfrak p}_i}\not\subset{\mathfrak P}_iS_{{\mathfrak p}_i}$. 
Donc pour 
$Q_i\in{\left({{\mathfrak p}_i}^{\delta_i}\right)}^{\beta_i}$ 
g\'en\'erique 
$Q_iS_{{\mathfrak p}_i}\not\subset{\mathfrak P}_iS_{{\mathfrak p}_i}$. 
Donc $({\mathfrak P}_iS_{{\mathfrak p}_i}:Q_i)$ est un id\'eal strict de
$S_{{\mathfrak p}_i}$. Donc $({\mathfrak P}_iS_{{\mathfrak p}_i}:Q_i)
\subset{\mathfrak p}_iS_{{\mathfrak p}_i}$.
Comme $({\mathfrak P}_i:Q_i)$ est primaire, nous en d\'eduisons 
$({\mathfrak P}_i:Q_i)\subset{\mathfrak p}_i$.
R\'eciproquement, comme $Q_i\in({{\mathfrak p}_i})^{\beta_i}$, 
nous avons $Q_i{\mathfrak p}_i\subset ({{\mathfrak p}_i})^{\beta_i+1}
\subset{\mathfrak P}_i$, donc ${\mathfrak p}_i\subset({\mathfrak P}_i:Q_i)$.
Finalement ${\mathfrak p}_i=({\mathfrak P}_i:Q_i)$. Posons $Q=\prod_iQ_i$. 
L'assertion~2 r\'esulte de ce que  
$I_W=\bigcap_i{\mathfrak p}_i=\bigcap_i({\mathfrak P}_i:Q_i)=(I_V:Q)$.  
\cucu

\begin{lemme}\label{bamu}
Il existe une
constante $C$ qui ne d\'epend que de 
$\varepsilon$, $S$, $u$, $s$ et $m$, telle que les id\'eaux
$I_W$ et ${\mathfrak q}$ sont engendr\'es en degr\'e
inf\'erieur  \`a $Cb^{a}$, o\`u $a=2^{A-1}$.    
\end{lemme}

\preuve
D'apr\`es l'assertion 2 du lemme~\ref{alh-}, nous avons
${\mathfrak I}_V=({\mathcal F}_0,\dots,{\mathcal F}_{A-n-1}:
{\mathcal G})$, o\`u les ${\mathcal F}_i$ sont en
intersection compl\`ete et v\'erifient 
$\deg {\mathcal F}_i\leq(1+\varepsilon)b$,
et o\`u $\deg {\mathcal G}\leq(A-n+1)(1+\varepsilon)b+C'$, 
o\`u $C'$ est la constante donn\'ee par le lemme~\ref{alh-}. 
Notons ${\mathfrak I}_W$ la pr\'einage de 
$I_W\subset S$ dans ${\mathfrak S}$.
L'assertion~2 du lemme~\ref{reduction} implique
${\mathfrak I}_W=({\mathcal F}_0,\dots,{\mathcal F}_{A-n-1}:
{\mathcal G}{\mathcal Q})$, o\`u ${\mathcal Q} $ 
d\'esigne une pr\'eimage de $Q\in S$ dans ${\mathfrak S}$, et o\`u
$\deg {\mathcal G}{\mathcal Q}=
\deg {\mathcal G}+\deg {\mathcal Q}\leq(A-n+1)(1+\varepsilon)b+
C'+(1+\varepsilon)b\leq (A-n+2)(1+\varepsilon)b+C'$. 
La suite exacte
\begin{eqnarray*}
0\longrightarrow {\mathfrak S}/{\mathfrak I}_W
[-\deg {\mathcal G}{\mathcal Q}]
   \stackrel{\cdot {\mathcal G}{\mathcal Q}}{\longrightarrow}
{\mathfrak S}/{({\mathcal F}_0,\dots,{\mathcal F}_{A-n-1})}
   \longrightarrow {\mathfrak S}/
{({\mathcal F}_0,\dots,{\mathcal F}_{A-n-1},{\mathcal G}{\mathcal Q})}
  \longrightarrow 0
\end{eqnarray*}
implique que la r\'egularit\'e au sens de Mumford-Castelnuovo
({\it cf.} \cite{bamu}, d\'efinition 3.2) de ${\mathfrak I}_W$
est major\'ee par celle des id\'eaux $({\mathcal F}_0,\dots,{\mathcal F}_{A-n-1})$ et
$({\mathcal F}_0,\dots,{\mathcal F}_{A-n-1},{\mathcal G}{\mathcal Q})$. Or ces id\'eaux sont engendr\'es en
degr\'e inf\'erieur  \`a 
$(A-n+2)(1+\varepsilon)b+C'$. D'apr\`es
\cite{bamu}, th\'eor\`eme 3.7, la r\'egularit\'e de 
${\mathfrak I}_W$ est donc major\'ee par  
$(2(A-n+2)(1+\varepsilon)b+C')^{2^{A-2}}$~; 
donc ${\mathfrak I}_W$ est engendr\'e en degr\'e inf\'erieur \`a   
$(2(A-n+2)(1+\varepsilon)b+C')^{2^{A-2}}$, ce qui
d\'emontre le lemme~\ref{bamu} pour $I_W$. 

Pour tout $i$ nous notons ${\mathcal G}_i$ une pr\'eimage dans 
${\mathfrak S}$
d'un \'el\'ement g\'en\'erique de ${\mathfrak p}_i^{\delta_i}$. 
Nous posons ${\mathcal Q}_i={\mathcal Q}\prod_{j\neq i}{\mathcal G}_j$~: 
nous avons alors 
$\deg {\mathcal Q}_i\leq\deg V$ et ${\mathfrak p}_i=(I_V:{\mathcal Q}_i)$. 
En rempla\c cant ${\mathcal Q}$ par ${\mathcal Q}_i$ dans le
raisonnement ci-dessus nous obtenons que l'id\'eal ${\mathfrak p}_i$  
est engendr\'e en degr\'e inf\'erieur  \`a   
$(2(A-n+2)(1+\varepsilon)b+C')^{2^{A-2}}$.
Comme $\sum_i\beta_i\leq(1+\varepsilon)b$, l'id\'eal
${\mathfrak q}=\prod_i({\mathfrak p}_i)^{\beta_i}$ 
est donc engendr\'e en degr\'e inf\'erieur  \`a
$(1+\varepsilon)b
(2(A-n+2)(1+\varepsilon)b+C')^{2^{A-2}}$, ce qui
d\'emontre le lemme~\ref{bamu} pour ${\mathfrak q}$.
\cucu

\section{Preuve de la proposition~\ref{I^2am}}\label{preuveI^2am}

Nous invitons le lecteur \`a supposer $V$ r\'eduit en premi\`ere 
lecture (nous avons alors ${\mathfrak q}=S$).

Soit $\gamma$ le plus petit entier tel que les id\'eaux $I_W$ 
et ${\mathfrak q}$ sont engendr\'es en degr\'e inf\'erieur ou 
\'egal \`a $\gamma$. D'apr\`es l'assertion~2 de la proposition~\ref{alh},
nous avons $\gamma\leq Cb^{a}$, o\`u $C$ est une constante 
qui ne d\'epend que de $Y$, de $\anneau_Y(1)$ et de $\varepsilon$.

Sans perdre en g\'en\'eralit\'e nous supposons $\varepsilon<\frac{1}{2N}$.
Il existe alors une constante $C\in\RR_+^*$ qui ne d\'epend que de $Y$, 
de $\anneau_Y(1)$ et de $\varepsilon$, telle que pour tout $d\geq Cb^a$
et pour tout $r\in\NN^*$ nous avons 
$\frac{1+\varepsilon}{\varepsilon} b\leq d$,
$2(r+1)\gamma\leq rd$ et
\begin{eqnarray} \label{ineqaffreuse}
   (n+r+1)d-m-2(r+1)\gamma > (n+1)(d+s)+(N-n+1)\varepsilon rd +C',
\end{eqnarray}
o\`u $C'$ est la constante d\'efinie \`a la section~\ref{subalge}.
Nous supposons $d\geq Cb^a$.

La preuve de la proposition~\ref{I^2am} proc\`ede par 
approximations successives~: le lemme~\ref{sous-lemme} ci-dessous
peut \^etre consid\'er\'e comme une version affaiblie de la 
proposition~\ref{I^2am}. 

\begin{lemme}\label{vieE_r} Pour tout $L\in S^{\leq 2(r+1)\gamma}$ soit 
$L\in E_r$, soit l'id\'eal $(E_r:L)$ est distinct de $S$ en 
degr\'e $(n+1)(d+s)+(N-n+1)\varepsilon rd +C'$.
\end{lemme}

\preuve
Supposons $L\not\in E_r$, donc $(E_r:L)\neq S$.
Alors il existe $i\in\{1,\dots,u\}$ tel que $(E_{r,i}:L)\neq S$~:
l'id\'eal $(E_{r,i}:L)$ est Gorenstein de degr\'e du socle
$(n+r+1)d-m-\deg L$. Comme $\deg\,L\leq 2(r+1)\gamma$, 
l'in\'egalit\'e~(\ref{ineqaffreuse}) implique 
que $(E_{r,i}:L)$, donc \`a fortiori $(E_{r}:L)$ sont 
distincts de $S$ en degr\'e $(n+1)(d+s)+(N-n+1)\varepsilon rd +C'.$
\cucu

\begin{lemme} \label{spb} Pour tout
$L\in E_{r-1}$ dont les d\'eriv\'ees partielles
$L_{\widetilde{v}}L$, $\widetilde{v}\in \widetilde{T}$ appartiennent 
\`a $E_{r-1}$, nous avons 
$l_j(E_r:L)\leq d+s$ pour tout $j\in\{0,\dots,N+1\}$.
\end{lemme}

\preuve
Comme $L\in E_{r-1}$, l'assertion~1 de la proposition~\ref{3prop} 
implique $F\in(E_r:L)$. Pour tout $\widetilde{v}\in \widetilde{T}$,
comme $L_{\widetilde{v}}L\in E_{r-1}$, l'assertion~3 de la 
proposition~\ref{3prop} implique alors $L_{\widetilde{v}}F\in(E_r:L)$.
Donc $J_F\subset(E_r:L)$ et d'apr\`es la section~\ref{jaco}, 
l'id\'eal $(E_r:L)$ est sans point base en degr\'e $d+s$. 
\cucu

\begin{lemme}\label{sous-lemme}
Pour $r\geq 1$ nous avons 
${\left(E_{r-1}^{\leq (r+1)\gamma}\right)}^2\subset E_{r}$.
\end{lemme}

\preuve
Soit $L\in (E_{r-1})^2$ tel que $\deg\,L\leq 2(r+1)\gamma$.
D'une part, l'assertion 2 de la proposition~\ref{3prop} implique 
$\Theta\subset(E_r:L)$, donc $\dim\left(S^d/(E_{r}:L)\right)\leq
\dim\left(S^d/\Theta\right)\leq b\frac{d^n}{n!}$. 
Comme $2(r+1)\gamma\leq rd$, l'id\'eal $(E_r:L)$ est intersection 
des id\'eaux $(E_{r,i}:L)$ de degr\'e du socle
$(n+r+1)d-\deg\,L\geq(n+1)d-m$ et l'id\'eal $(E_r:L)$ satisfait les 
hypoth\`eses~2 et~3 de la donn\'ee alg\'ebrique~\ref{doal}. D'apr\`es le 
lemme~\ref{llh} nous avons donc
$l_j(E_r:L)\leq\varepsilon d\leq\varepsilon rd$ pour tout 
$j\in\{0,\dots,N-n\}$.
D'autre part, le lemme~\ref{spb} implique
$l_j(E_r:L)\leq d+s$ pour tout $j\in\{N-n+1,\dots,N+1\}$.
L'id\'eal $(E_{r}:L)$ contient donc une intersection compl\`ete 
sans point base dont la somme des degr\'es est inf\'erieure  \`a 
$(N-n+1)\varepsilon rd+(n+1)(d+s)$.
La section~\ref{subalge} et le lemme~\ref{vieE_r} impliquent alors 
$L\in E_r$.
\cucu

\bigskip

\noindent
{\bf Fin de la preuve de la proposition~\ref{I^2am}.}
Nous proc\'edons par r\'ecurrence sur $r$.

Pour $r=0$ la proposition~\ref{I^2am} r\'esulte de la 
proposition~\ref{alh} et de ce que ${\mathfrak q}I_W\subset I_V$.

Supposons $r>0$ et l'\'enonc\'e vrai pour $r-1$~: nous avons
${\mathfrak q}(I_W)^{r+1}\subset E_{r-1}$, et comme 
l'id\'eal ${\mathfrak q}(I_W)^{r}$
est engendr\'e en degr\'e inf\'erieur \`a 
$(r+1)\gamma$, ${\mathfrak q}(I_W)^{r+1}\subset E_{r-1}^{\leq(r+1)\gamma}$.

Soit $L\in{\mathfrak q}(I_W)^{r+1}$~; nous pouvons supposer 
$\deg L\leq(r+2)\gamma$, donc $\deg L\leq 2(r+1)\gamma$. 
D'une part l'hypoth\`ese de r\'ecurrence et le lemme~\ref{sous-lemme}
donnent
${\left({\mathfrak q}\right)}^2(I_W)^{2r}\subset 
{\left(E_{r-1}^{\leq (r+1)\gamma}\right)}^2\subset 
{E_{r}}$, et donc
${\mathfrak q}(I_W)^{r-1}\subset (E_{r}:L)$. L'assertion~1 du
lemme~\ref{reduction} et le fait que $W$ est de 
dimension pure~$n$ et de degr\'e inf\'erieur  
\`a $(1+\varepsilon)b$ impliquent alors
$l_j(E_{r}:L)\leq r(1+\varepsilon)b\leq\varepsilon rd$
pour tout $j\in\{0,\dots,N-n\}$. 
D'autre part $L$ et ses d\'eriv\'ees partielles 
$L_{\widetilde v}L$, ${\widetilde v}\in{\widetilde T}$,
appartiennent \`a ${\mathfrak q}(I_W)^{r}$, donc \`a $E_{r-1}$. 
Le lemme~\ref{spb} implique alors $l_j(E_r:L)\leq d+s$ pour 
tout $j\in\{N-n+1,\dots,N+1\}$.
L'id\'eal $(E_{r}:L)$ contient donc une intersection compl\`ete 
sans point base dont la somme des degr\'es est inf\'erieure  \`a 
$(N-n+1)\varepsilon rd+(n+1)(d+s)$.
La section~\ref{subalge} et le lemme~\ref{vieE_r} impliquent alors 
$L\in E_r$.
\cucu

\section{Preuve de la proposition~\ref{Zred}}\label{preuveZred}

\subsection{Pr\'eliminaires}\label{preZ}
Les notation sont celles de la section~\ref{introNL}.
Le faisceau inversible $\anneau_Y(1)$ d\'efinit un plongement
$Y\hookrightarrow\PP^{A-1}_{\C}$, o\`u $A=\hh^0(Y,\anneau_Y(1))$. 
Soit ${\bf G}={\rm Gr}_{\C}\,(A-n-2,A)$ la grassmanienne~; pour tout
$l\in{\bf G}$, soit $L\subset\PP^{A-1}_{\C}$ l'espace lin\'eaire de
codimension $n+2$ associ\'e et pour tout sch\'ema
$Z\subset\PP^{A-1}_{\C}\setminus L$, soit 
$C(Z,l)\subset\PP^{A-1}_{\C}$ le c\^one projectif de sommet $L$ et de
base $Z$.   

\begin{lemme} \label{avantchamp}
Soit $M\subset X_F$ une vari\'et\'e de dimension $n-1$. 
Alors pour $l\in{\bf G}$ g\'en\'erique, les vari\'et\'es 
$C(M,l)$ et $X_F$ sont transverses en tout point de $M$.
\end{lemme}

\preuve 
Posons
\begin{eqnarray*}
\Gamma_1&=&\{(x,l)\in M\times {\bf G}\ |\ x\in L\}\qquad{\rm et}\\
\Gamma_2&=&\{(x,l)\in M\times {\bf G}\setminus\Gamma_1\ |\ 
T_xC(x,l)+T_xX_F\neq T_x\PP_{\C}^{A-1}\}.
\end{eqnarray*}

D'une part nous avons 
$\dim\Gamma_1=\dim{\bf G}+\dim M-\codim(L,\PP_{\C}^{A-1})=\dim{\bf G}-3$.

D'autre part, pour tout $(x,l)\in X_F\times {\bf G}\setminus\Gamma_1$,
le sous-espace vectoriel $T_xC(x,l)+T_xX_F$ de $T_x\PP_{\C}^{A-1}$ est
d\'efini par une matrice de $A-1$ lignes et $(A-n-2)+N$ colonnes dont les
coefficients d\'ependent de $(x,l)$~; la condition
$T_xC(x,l)+T_xX_F\neq T_x\PP_{\C}^{A-1}$ \'equivaut \`a l'annulation de
ses mineurs de taille $(A-1)\times(A-1)$, ce qui donne $N-n$
conditions ind\'ependantes. Nous avons donc
$\dim\Gamma_2=\dim{\bf G}+\dim M-(N-n)=\dim{\bf G}-(N-2n+1)$.

Posons $\Gamma=\Gamma_1\cup\Gamma_2$. Comme $N\geq 2n$, nous avons
donc $\dim\Gamma<\dim{\bf G}$~; la projection de $\Gamma$ sur ${\bf G}$
a donc pour image un ferm\'e strict de ${\bf G}$ et le
lemme~\ref{avantchamp} est vrai pour tout $l\in {\bf G}$ n'appartenant pas
\`a ce ferm\'e.
\cucu

\medskip

Nous adoptons les notations de la section~\ref{derivation}. Pour tout
$y\in Y$, $\widetilde{y}\in\pi^{-1}(y)$ et 
$\widetilde{v}\in\widetilde{T}$, le vecteur
$\pi_*\widetilde{v}(\widetilde{y})\in T_yY$ ne d\'epend du choix de 
$\widetilde{y}$ que par la multiplication par un scalaire. 
Nous pouvons donc poser
$$\widetilde{T}(l)=\{\widetilde{v}\in\widetilde{T}\ |\
\forall y\in Y,\ \pi_*\widetilde{v}(\widetilde{y})\in T_y(C(y,l)\cap Y)\}.$$
C'est un $S$-module de type fini~; il est donc engendr\'e
en degr\'e inf\'erieur  \`a un entier relatif $t(l)$. Notons $t$ la
valeur de $t(l)$ pour $l\in {\bf G}$ g\'en\'eral~; c'est un entier qui ne
d\'epend que de $Y$ et de $\anneau(1)$.

\begin{lemme} \label{champ}
Soit $Z$ une sous-vari\'et\'e irr\'eductible de $Y$ de dimension pure
$n$ et de degr\'e $\alpha$. Alors si $Z\not\subset X_F$, 
pour $l\in{\bf G}$ g\'en\'eral, l'id\'eal
$\left(I_Z,F,L_{\widetilde{v}}F\ |\ 
\widetilde{v}\in\widetilde{T}(l)\right)$
co\"{\i}ncide avec $S$ en degr\'e
$(N-n+1)\alpha+d+n(d+t)+C'$, o\`u $C'$ est la constante d\'efinie \`a la 
section~\ref{subalge}.
\end{lemme}

\preuve
Posons $I=\left(I_Z,F,L_{\widetilde{v}}F\ |\ 
\widetilde{v}\in\widetilde{T}(l)\right)$.
Comme $Z$ est de dimension pure $n$, pour $l\in {\bf G}$ g\'en\'erique la
vari\'et\'e $C(Z,l)\cap Y$ est une hypersurface de $Y$ de degr\'e
$\alpha$ dont l'\'equation appartient \`a $I_Z$~; nous en d\'eduisons
$l_i(I)\leq l_i(I_Z)\leq\alpha$ pour $i\in\{0,\dots,N-n\}$.

Comme $Z$ est irr\'eductible et que $Z\not\subset X_F$, nous avons
$\dim Z\cap X_F=n-1$, et comme $F\in I$, nous en d\'eduisons 
$l_{N-n+1}\leq d$. 

Enfin, le lemme~\ref{avantchamp} appliqu\'e \`a
$M=Z\cap X_F$ implique que si $l\in {\bf G}$ est g\'en\'erique 
alors pour tout $y\in Z\cap X_F$ nous avons 
$T_yX_F\not\subset T_y(C(y,l)\cap Y)$. 
Donc si $\widetilde{T}(l)$ est engendr\'e en degr\'e inf\'erieur 
\`a $t$, il existe
$\widetilde{v}\in\widetilde{T}(l)^{\leq t}$ et 
$\widetilde{y}\in\pi^{-1}(y)$ tels que 
$\pi_*\widetilde{v}(\widetilde{y})\not\in T_yX_F$ et donc
$L_{\widetilde{v}}F(y)\neq 0$. 
Donc pour $l\in{\bf G}$ g\'en\'eral l'id\'eal
$\left(L_{\widetilde{v}}F\ |\ 
\widetilde{v}\in\widetilde{T}(l)^{\leq t}\right)$
ne s'annule pas sur $Z\cap X_F$ et nous avons 
$l_i(I)\leq d+t$ pour $i\in\{N-n+2,\dots,N+1\}$.

Finalement $\sum_{i=0}^{N+1}l_i\leq
(N-n+1)\alpha+d+n(d+t)$, ce qui implique le lemme~\ref{champ}.
\cucu

\subsection{Preuve de la proposition~\ref{Zred}}
Nous adoptons les notations des sections~\ref{alg},~\ref{sebamu} 
et~\ref{preZ}. 

Il existe une constante $C\in\RR_+^*$ qui ne d\'epend que de $Y$, 
de $\anneau_Y(1)$ et de $\varepsilon$, qui v\'erifie
la proposition~\ref{I^2am} et telle que pour tout 
$d\geq Cb^a$ nous avons 
\begin{eqnarray} 
(N-n+1)\alpha+d+n(d+t)+C'&<&(n+2)d-m-r\label{ineq1},\quad {\rm et}\\
r+t&\leq& d-(A-n-1)(1+\varepsilon)b+C''+t,
\label{ineq2}
\end{eqnarray}
o\`u $C'$ est la constante d\'efinie \`a la 
section~\ref{subalge}, o\`u $C''$ est la constante du lemme~\ref{alh-} et
o\`u $r\in \NN^*$ v\'erifie $r\leq 2(1+\varepsilon)b$.
Nous supposons $d\geq Cb^a$.

Soit $Z\subset W$ une composante irr\'eductible~; nous 
devons montrer $Z\subset X_F$. Ceci r\'esulte du lemme~\ref{champ} 
ci-dessus, des lemmes~\ref{contenu} et~\ref{vie} ci-dessous 
et de l'in\'egalit\'e~(\ref{ineq1}). 
\cucu

\bigskip

Soit $Z'$ l'adh\'erence de $W\setminus Z$ dans $W$~: nous avons
$W=Z\cup Z'$.  Comme les sch\'emas $Z$, $Z'$ et $V({\mathfrak q})$ 
sont de dimension pure $n$, pour $l\in {\bf G}$ g\'en\'erique les 
sch\'emas $C(Z,l)$, $C(Z',l)$ et $C(V({\mathfrak q}),l)$ sont des 
hypersurfaces de $\PP^{A-1}_{\C}$ transverses \`a $Y$.
Posons $\alpha=\deg Z$ et $\alpha'=\deg Z'$. Alors
$\deg\,C(Z,l)=\alpha$, $\deg\,C(Z',l)=\alpha'$ et 
$\deg\,C(V({\mathfrak q}),l)=\beta$.
Notons $P_l\in S^{\alpha}$, $P'_l\in S^{\alpha'}$ et 
$Q_l\in S^{\beta}$ les \'equations respectives de leurs 
intersections avec $Y$ (d\'efinies \`a la multiplication par un 
scalaire pr\`es). Posons $R_l=P_l{P'_l}^2Q_l$ et $r=\deg\,R_l$. 
Nous avons $r=\alpha+2\alpha'+\beta\leq 2(1+\varepsilon)b$.

Le lecteur pourra supposer $V$ irr\'eductible et r\'eduit
en premi\`ere lecture~: si V est irr\'eductible alors $Z=W$ et
l'hypersurface $C(Z',l)$ est vide~; si $V$ est r\'eduit alors $W=V$ et 
l'hypersurface $C(V({\mathfrak q}),l)\cap Y$ est vide~;
donc si $V$ est irr\'eductible et r\'eduit alors
$P'_l=Q_l=1$, $r=\deg\,V$ et $R_l$ est une \'equation de $C(V,l)$.

\begin{lemme} \label{contenu} 
Pour $l\in{\bf G}$ g\'en\'erique $\left(I_Z,F,L_{\widetilde{v}}F\ |\ 
\widetilde{v}\in\widetilde{T}(l)\right)\subset(E_1:R_l)$.
\end{lemme}

\preuve
D'une part, pour tout $K\in I_Z$ nous avons $KP'_l\in I_W$, donc 
$KR_l\in {\mathfrak q}(I_W)^2$, donc $KR_l\in E_1$ par la 
proposition~\ref{I^2am}, donc $I_Z\subset(E_1:R_l)$. 
D'autre part, comme $P_lP'_lQ_l\in I_V$, donc $R_l\in I_V$, 
la proposition~\ref{alh} implique $R_l\in E_0$
et l'assertion~1 de la proposition~\ref{3prop} implique
$F\in(E_1:R_l)$. Enfin, comme $R_l$ s'annule sur $C(V,l)$,
pour tout $\widetilde{v}\in\widetilde{T}(l)$, 
$L_{\widetilde{v}}R_l$ s'annule sur $C(V,l)$,
donc $L_{\widetilde{v}}R_l\in I_V$.
La proposition~\ref{alh} implique $L_{\widetilde{v}}R_l\in E_0$ 
et l'assertion~3 de la proposition~\ref{3prop} implique 
$L_{\widetilde{v}}F\in(E_1:R_l)$.
\cucu

\begin{lemme} \label{vie} 
Pour $l\in{\bf G}$ g\'en\'eral l'id\'eal $(E_1:R_l)$ est distinct de 
$S$ en degr\'e $(n+2)d-m-r$.
\end{lemme}

\preuve
Nous avons $I_Z={\mathfrak p}_i$ pour un certain $i$. 
Alors $Q_l\not\in({\mathfrak p}_i)^{\beta_i+1}$, donc 
$R_l\not\in({\mathfrak p}_i)^{\beta_i+2}$, donc
il existe $\widetilde{v}\in\widetilde{T}^{\leq t}$ tel que 
$L_{\widetilde{v}}R_l\not\in({\mathfrak p}_i)^{\beta_i+1}$, 
donc $L_{\widetilde{v}}R_l$ ne s'annule pas sur 
$C(V(({\mathfrak p}_i)^{\beta_i+1}),l)$.
Comme pour $l\in{\bf G}$ g\'en\'erique
$C(V(({\mathfrak p}_i)^{\beta_i+1}),l)$
est une composante irr\'eductible de $C(V,l)$, 
$L_{\widetilde{v}}R_l$ ne s'annule pas sur $C(V,l)$ et 
$L_{\widetilde{v}}R_l\not\in I_V$. Comme 
$\deg\,L_{\widetilde{v}}R_l\leq r+t$, 
l'in\'egalit\'e~(\ref{ineq2}) et l'assertion~3 du lemme~\ref{alh-}
impliquent $L_{\widetilde{v}}R_l\not\in E_0$. 
Les assertions~1 et~3 de la proposition~\ref{3prop} impliquent
alors $R_l\not\in E_1$. Comme 
$E_1=\bigcap_{i=1}^u E_{1,i}$, il existe 
$i\in\{1,\dots,u\}$ tel que $R_l\not\in E_{1,i}$~:
l'id\'eal $(E_{1,i}:R_l)$ est Gorenstein de degr\'e du socle 
$(n+2)d-m-r$. Le lemme~\ref{vie} r\'esulte alors de ce que
$(E_1:R_l)\subset (E_{1,i}:R_l)$. 
\cucu

\medskip

\noindent
Ania Otwinowska\\
Laboratoire de Math\'ematiques,\
Universit\'e Paris-Sud - B\^at 425,\
91405 Orsay Cedex,\
France\\
e-mail~: Ania.Otwinowska@math.u-psud.fr

\begin{thebibliography}{999999999}
\bibitem[Ba-Mu]{bamu}
           {\sc D. Bayer, D. Mumford.}
           {\it What can be computed in algebraic geometry? }
           {Computational algebraic geometry and commutative algebra
           (Cortona 1991), 1-48, Sympos. Math. XXXIV, Cambridge
           Univ. Press, Cambridge (\oldstylenums{1996}).}
\bibitem[CCGH]{ivhs}
           {\sc J. Carlson, M. Green, Ph. Griffiths, J. Harris}
           {\it Infinitesimal Variations of Hodge Structures I.}
           {Comp. Math 50, p. 105-205 (\oldstylenums{1983}).}    
\bibitem[CDK]{CDK}
           {\sc E. Cattani, P. Deligne, A. Kaplan.}
           {\it On the locus of Hodge classes.}
           {Journal of the AMS, Vol 8, n.2 (\oldstylenums{1995}).}    
\bibitem[Gre]{NLgreen} 
           {\sc M. Green.}
           {\it Components of maximal dimension in the
           Noether-Lefschetz locus.}
           {J. Differential Geometry 29, p. 295-302 (\oldstylenums{1989}).}
\bibitem[Gri]{grI}
           {\sc Ph. Griffiths.}
           {\it On the periods of certain rational integrals I.}   
           {Annals of Math. (2) 90, p. 460-495
           (\oldstylenums{1969}).}
\bibitem[GH]{ivhs2}
           {\sc Ph. Griffiths, J. Harris}
           {\it Infinitesimal Variations of Hodge Structures II.}
           {Comp. Math 50, p. 207-265 (\oldstylenums{1983}).}    
\bibitem[O1]{NL} 
           {\sc A. Otwinowska.}
           {\it Composantes de petite codimension du lieu de 
           Noether-Lefschetz~: un argument asymptotique en faveur de
           la conjecture de Hodge pour les hypersurfaces.} 
           {J. Algebraic Geom. 12, p. 307-320 (\oldstylenums{2003}).}
\bibitem[O2]{topo} 
           {\sc A. Otwinowska.}
           {\it Monodromie d'une famille d'hypersurfaces contenant 
           un sous-sch\'ema fix\'e.} 
           {Soumis, (\oldstylenums{2000}).}
\bibitem[V1]{NLclaire}
           {\sc C. Voisin.}
           {\it Une pr\'ecision concernant le th\'eor\`eme de Noether.}
           {Math. Ann. 280, p. 605-611 (\oldstylenums{1989}).}
\bibitem[V2]{NLclaire2}
           {\sc C. Voisin.}
           {\it Composantes de petite codimension du lieu de 
           Noether-Lefschetz.} 
           {Comment. Math. Helvetici 64, p. 515-526 (\oldstylenums{1989}).}
\end{thebibliography}
\end{document}